\documentclass{article}

\newif\ifarxiv
\arxivtrue

\usepackage{graphicx}
\usepackage{amsmath,amssymb,amsthm,bm}
\usepackage[a4paper]{geometry}
\usepackage{algorithm,algorithmic}
\usepackage[colorlinks=false,pdfborder={0 0 0}]{hyperref}
\usepackage{authblk}
\usepackage{color}
\usepackage{multirow}

\graphicspath{{fig/}}
\DeclareGraphicsExtensions{.mps,.pdf,.eps,.png}

\newenvironment{keywords}{\medskip\textbf{Keywords:}}{}
\newenvironment{AMS}{\medskip\textbf{AMS subject classifications (2020).}}{}

\newtheorem{theorem}{Theorem}

\newtheorem{remark}{Remark}
\newtheorem{lemma}{Lemma}

\theoremstyle{plain}

\DeclareMathOperator*{\argmin}{\arg\min}

\DeclareMathOperator{\Span}{span}

\newcommand{\fro}{\mathsf F}

\newcommand*{\trans}{^{\top}}

\newcommand{\bmat}[1]{\begin{bmatrix}#1\end{bmatrix}}

\newcommand*{\norm}[1]{\Vert#1\rVert}
\newcommand*{\normbig}[1]{\big\Vert#1\big\rVert}

\newcommand*{\normF}[1]{\bigl\Vert#1\bigr\rVert_{\fro}}

\newcommand{\sigmin}{\sigma_{\min}}

\def\adots{\mathinner{\mkern2mu\raise1pt\hbox{.}\mkern2mu
    \raise4pt\hbox{.}\mkern2mu\raise7pt\hbox{.}\mkern1mu}}

\newcommand*{\macheps}{\bm u}

\newcommand*{\bigO}{O}

\newcommand*{\epsapplyS}{\varepsilon_{S}}
\newcommand*{\epsb}{\varepsilon_b}

\newcommand*{\epskappaW}{\omega_{\kappa(W)}}
\newcommand*{\epskapparW}{\omega_{\kappa(r, W)}}

\newcommand*{\epssls}{\varepsilon_{sls}}
\newcommand*{\epsW}{\varepsilon_{W}}
\newcommand*{\epsx}{\varepsilon_{x}}

\newcommand*{\krylov}{\mathcal{K}}

\title{On the numerical stability of sketched GMRES}
\author[1]{Liam Burke}
\author[1]{Erin Carson}
\author[1]{Yuxin Ma}
\affil[1]{Department of Numerical Mathematics, Faculty of Mathematics and Physics, Charles University, Sokolovsk\'{a} 49/83, 186 75 Praha 8, Czechia

\texttt{Email: liam.burke@matfyz.cuni.cz, carson@karlin.mff.cuni.cz, yuxin.ma@matfyz.cuni.cz}}

\begin{document}
\maketitle

\begin{abstract}
We perform a backward stability analysis of preconditioned sketched GMRES [Nakatsukasa and Tropp, SIAM J. Matrix Anal. Appl, 2024] for solving linear systems $Ax=b$, and show that the backward stability at iteration $i$ depends on the conditioning of the Krylov basis $B_{1:i}$ as long as the condition number of $A B_{1:i}$ can be bounded by \(1/\bigO(\macheps)\), where \(\macheps\) is the unit roundoff. Under this condition, we demonstrate that the stability of the sketched GMRES is mainly affected by the condition number of $B_{1:i}$, which can be mitigated by using the restarting technique.
We also derive sharper bounds that
better capture the attainable backward error especially for cases when the basis $B_{1:i}$ is very ill-conditioned, which has been observed in the literature but
not yet explained theoretically.
We present numerical experiments to demonstrate the conclusions of our analysis, and also show that adaptively restarting where appropriate allows us to recover backward stability in sketched GMRES. 


\end{abstract}

\section{Introduction}
\label{sec:introduction}
In this paper, we analyze the numerical stability of sketched GMRES (sGMRES) \cite{NT2024} for solving nonsymmetric and nonsingular linear systems of the form  
\begin{equation}\label{eq:linear_system}
    A x = b, \enspace A \in \mathbb{R}^{n \times n}
\end{equation}
with a given matrix \(A\) and a right-hand side \(b\).

Although standard GMRES \cite{SaadSchultz1986} is one of the most popular methods to solve \eqref{eq:linear_system}, it comes at the expense of constructing a matrix $V_{1:i} \in \mathbb{R}^{n \times i}$ with columns that form an orthonormal basis for the Krylov subspace 
\[\krylov_i(A, r^{(0)}) = \Span\{r^{(0)}, A r^{(0)}, \dotsc, A^{i-1} r^{(0)}\},
\]
where $i$ is the iteration number, and \(r^{(0)} = b - A x^{(0)}\) is the initial residual corresponding to an initial approximation $x^{(0)}$.

Sketched GMRES was proposed in \cite{NT2024} as an attractive alternative to GMRES that exploits the technique of randomized sketching \cite{martinsson2020randomized} to avoid the high cost of constructing an orthonormal basis. Instead, the  basis is only partially orthogonalized, thereby reducing the arithmetic cost from \(\bigO(ni^2)\) to \(\bigO(i^3+ni\log i)\).

It has been observed in numerical experiments in~\cite{NT2024} that sGMRES can fail once the condition number of the basis $B_{1:i}$ becomes too large, and methods such as the \textit{sketch-and-select Arnoldi} algorithm \cite{guttel2024sketch} have been proposed as a means to ensure the condition number of the basis $B_{1:i}$ does not grow too rapidly. Indeed, numerical experiments in \cite{guttel2024sketch} showed that working with a well-conditioned basis improved the overall convergence behavior of sketched GMRES.

In this paper, we rigorously analyze the backward stability of preconditioned sketched GMRES, and justify the experimental observations made in the literature \cite{NT2024, guttel2024sketch}.
Namely, we show that the backward error of sketched GMRES depends on \(\kappa(Z_{1:i})\) as long as \(\kappa(AZ_{1:i})\) is not too large, where \(Z_{1:i} = M_R^{-1}B_{1:i}\) with the right preconditioner \(M_R\).
We then provide a backward stability result for a restarted version of sketched GMRES, which under certain assumptions, allows us to recover the numerical stability of sketched GMRES when \( \kappa(Z_{1:i})\) grows large, suggesting we should monitor the conditioning of the basis throughout the iterations and restart when necessary. 
We present numerical experiments to demonstrate the numerical behavior we expect to observe in sketched GMRES as predicted by our new analysis, and demonstrate the improvements in the backward error that can be obtained from the adaptive restart approach.

We briefly outline the notation and some basic definitions used throughout this paper. We use MATLAB indexing to denote submatrices. For example, we use $X_{1:i}$ to denote the first $i$ columns of $X$. We use superscripts to denote iteration numbers. For example, $x^{(i)}$ denotes the approximate solution of $A x = b$ at the $i$-th iteration of GMRES or sketched GMRES. We use $\norm{\cdot}$ to denote the vector $2$ norm and $\normF{\cdot}$ to denote the Frobenius norm. We use \(\macheps\) to denote the unit roundoff. We also use $\hat{\cdot}$ to denote computed quantities.
For example, $\hat{x}$ denotes the computed value of the vector $x$ in finite precision.
The (normwise) backward error of $\hat{x}$ is defined as
\[
\min \{ \epsilon : (A + \Delta A) \hat{x} = b + \Delta b, \normF{\Delta A} \leq \epsilon \normF{A}, \norm{\Delta b} \leq \epsilon \norm{b} \} = \frac{\norm{b - A \hat{x}}}{\normF{A} \norm{\hat{x}} + \norm{b}}.
\]

The remainder of this paper is organized as follows.
In Section~\ref{sec:GMRES}, we give an overview of GMRES and the MOD-GMRES framework \cite{BHMV2024} for performing a backward stability analysis of different GMRES variants. We then discuss preconditioned sketched GMRES in Section \ref{sec:sketched_GMRES}, and analyze the backward stability of preconditioned sketched GMRES and restarted sketched GMRES, respectively, in Sections~\ref{sec:stability} and~\ref{sec:restart-sGMRES}.
Numerical experiments are presented in Section~\ref{sec:experiments}.

\section{Preconditioned GMRES}
\label{sec:GMRES}
In this section, we briefly outline the basics of preconditioned GMRES with a left preconditioner \(M_L\) for solving linear systems of the form~\eqref{eq:linear_system}. Beginning with an initial solution approximation $x^{(0)}$ and associated initial residual $r^{(0)} = b - A x^{(0)}$, preconditioned GMRES iteratively generates a matrix $V_{1:i} \in \mathbb{R}^{n \times i}$ with columns that form an orthonormal basis for the Krylov subspace \(  \krylov_i(M_L^{-1} A, M_L^{-1} r^{(0)}) \).
The orthonormal basis $V_{1:i}$ is typically constructed via the Arnoldi process, and satisfies the Arnoldi relation
\begin{align}
     M_L^{-1} A V_{1:i} = V_{1:(i+1)}\underline{H_i},
    \label{eq:Arnoldi}
\end{align}
    where $\underline{H_i}\in\mathbb{R}^{(i+1)\times i}$ is an upper Hessenberg matrix that stores the coefficients of the orthogonalization. The orthonormality of $V_{1:i}$ then allows for the GMRES least squares solution,
\begin{align}\label{eq:gmres_min}
        t^{(i)} 
        = 
        \argmin_{ t\in\mathcal{K}_i( M_L^{-1} A, M_L^{-1} r^{(0)})}
        \norm{M_L^{-1} r^{(0)} -  M_L^{-1} A t},
\end{align}
to be computed by solving a small  $(i+1) \times i$  least squares problem,
    \begin{align*}
        y^{(i)}
        =
        \argmin_{ y \in \mathbb{R}^i }
        \norm{M_L^{-1} r^{(0)} -  M_L^{-1} A V_{1:i} y^{(i)}}
        = 
        \argmin_{y \in \mathbb{R}^i}
        \norm{\| M_L^{-1} r^{(0)} \| e_1 - \underline{H_{1:i}} y},
\end{align*}
and setting $t^{(i)} = V_{1:i} y^{(i)}$. The solution approximation is then updated via $x^{(i)} = x^{(0)} + t^{(i)}$.
If a right preconditioner \(M_R\) is also taken into account, it can be integrated with the basis \(Z_{1:i}\) such that \(Z_{1:i} = M_R^{-1}V_{1:i}\).

There are many algorithmic variants of GMRES arising from different implementation options, such as the choice of preconditioner (if any), or the choice of orthogonalization method used within the Arnoldi procedure. Furthermore, for optimal performance it is common to exploit techniques such as deflation, mixed precision, or randomization within GMRES. In order to understand the finite precision behavior of any given GMRES variant, a backward error analysis framework was developed in \cite{BHMV2024}, which simplifies the process of deriving bounds for the attainable normwise backward and forward errors of the computed solutions of different GMRES variants. The framework is built from a backward error analysis of a generic GMRES algorithm known as \textit{modular GMRES} or \textit{MOD-GMRES}, that is composed of three elemental operations:
\begin{enumerate}
    \item Generate the preconditioned Krylov basis \(Z_{1:i} = M_R^{-1} V_{1:i}\) with the Krylov basis \(V_{1:i}\), and compute \(W_{1:i} = M_L^{-1} A Z_{1:i}\). \label{eq:mod-line-C}
    \item Solve the least squares problem $ y^{(i)} = \argmin_{y}\norm{\tilde b - W_{1:i} y}$ with \(\tilde b = M_L^{-1} b\).
    \item Compute the approximate solution \(x^{(i)} = Z_{1:i} y^{(i)}\). \label{eq:mod-line-x}
\end{enumerate}
Without loss of generality, we assume \(x^{(0)} = 0\) in this framework and our analysis of sGMRES in Section~\ref{sec:stability}.
By specializing these operations and meeting their assumptions, as outlined in the framework, modular GMRES can describe many GMRES implementations and variants, and can be used to perform their backward error analysis.
For example, well-known results from \cite{drkovsova1995numerical} were derived in \cite{BHMV2024} using the framework. In particular, it was shown that GMRES with Householder orthogonalization (HH-GMRES) produces a computed solution whose
backward error is on the order of the unit roundoff, and is thus backward stable. 
Similarly, the backward stability result for GMRES with modified Gram-Schmidt (MGS-GMRES), derived in \cite{paige2006modified}, was also recovered using the framework in \cite{BHMV2024}.

\section{Preconditioned Sketched GMRES}\label{sec:sketched_GMRES}
Before introducing sketched GMRES, we first recall the basic definition of the random sketching technique.
For $\epsilon \in (0,1)$, a sketching matrix \(S\in\mathbb{R}^{s\times n}\) is an \textit{$\epsilon$-subspace embedding} for the subspace \(\mathcal{V}\subseteq \mathbb{R}^{n\times m}\) if \(S\) satisfies that
\begin{equation}
(1-\epsilon) \| v \|^2 \leq \| S v \|^2 \leq (1+\epsilon) \|v\|^2,
\label{eq:sketch}
\end{equation}
for any $v \in \mathcal{V}$.
In practice, the sketching matrix \(S\) may not be explicitly available and needs to be drawn at random to achieve~\eqref{eq:sketch} with high probability.

Sketched GMRES aims to avoid the \(\bigO(n i^2)\) arithmetic cost of performing an expensive full orthogonalization of the basis vectors in standard GMRES. Using the technique of randomized sketching, the GMRES minimization problem \eqref{eq:gmres_min}
is replaced by the sketched problem
\begin{equation}\label{eq:sgmres_min}
y^{(i)} = \argmin_{y \in \mathbb{R}^{i}}\norm{S M_L^{-1} r^{(0)} - S M_L^{-1} A Z_{1:i} y},
\end{equation}
for a preconditioned basis \(Z_{1:i} = M_R^{-1} B_{1:i}\) with a (non-orthogonal) Krylov basis $B_{1:i}$ and sketching operator $S\in\mathbb{R}^{s\times n}$ with $i < s\ll n$.
The sketched problem \eqref{eq:sgmres_min} can then be constructed and solved inexpensively, even without orthonormality of the basis \( B_{1:i}\). The solution of the sketched problem serves as a good approximation of the original problem as long as the matrix $S$ is chosen to be an $\epsilon$-subspace embedding for the Krylov subspace $\mathcal{K}_{i}(M_L^{-1} A M_R^{-1}, M_L^{-1} r^{(0)})$, i.e., $S$ satisfies~\eqref{eq:sketch}
for $\epsilon \in (0,1)$, and every $v \in \mathcal{K}_{i}(M_L^{-1} A M_R^{-1}, M_L^{-1} r^{(0)})$.
Additionally, \(S\) must be constructed without prior knowledge of the Krylov subspace $\mathcal{K}_{i}(M_L^{-1} A M_R^{-1}, M_L^{-1} r^{(0)})$ and should typically be a valid embedding for a broad subspace, specifically \(\mathbb{R}^{n\times i}\).
In practice, the sketched least squares problem \eqref{eq:sgmres_min} is computed by solving a triangular linear system 
\begin{equation}\label{eq:triangular_system}
T_{1:i} y^{(i)} = U_{1:i}\trans (S M_L^{-1} r^{(0)}),
\end{equation}
where the orthonormal matrix \(U_{1:i}\) and the triangular matrix \(T_{1:i}\) come from the QR factorization of the reduced matrix \(C_{1:i} = S M_L^{-1} A Z_{1:i}\), i.e., \(C_{1:i} = U_{1:i} T_{1:i}\).

A detailed outline of one cycle of preconditioned sketched GMRES is presented in Algorithm~\ref{alg:sGMRES}, and a restarted version is given in Algorithm~\ref{alg:restarted_sGMRES}.

\begin{algorithm}[H]
\begin{algorithmic}[1]
    \caption{Preconditioned sketched GMRES (sGMRES) \label{alg:sGMRES}}
    \REQUIRE
    A matrix \(A \in \mathbb R^{n\times n}\), a right-hand side \(b \in \mathbb R^{n}\), an initial approximation \(x^{(0)} \in \mathbb R^{n}\), the maximal number of iterations \(m\), a subspace embedding \(S \in \mathbb R^{s\times n}\) (\(s>m\)), left and right preconditioners \(M_L\) and \(M_R\).
    \ENSURE
    A computed solution \(x \in \mathbb R^{n}\) approximating the solution of \(A x = b\).
    \STATE \(r^{(0)}\gets b - Ax^{(0)}\) and \(g \gets S M_L^{-1} r^{(0)}\). \label{algo-line:computer}
    \STATE \(B_{1} \gets r^{(0)}/\norm{r^{(0)}}\).
    \FOR{\(i = 1:m\)}
        \STATE Generate the \(i\)-th column \(B_{i}\) of the basis by Arnoldi process, and \(Z_{i}\gets M_R^{-1} B_i\).\STATE \(W_{i} \gets M_L^{-1} A Z_{i}\). \label{algo-line:computeW}
        \STATE Sketch reduced matrix \(C_{1:i} \gets S W_{1:i}\). \label{algo-line:computeC}
        \STATE Compute the QR factorization \(C_{1:i} = U_{1:i} T_{1:i}\) using \(C_{1:i-1} = U_{1:i-1} T_{1:i-1}\), where \(U_{1:i}\) is an orthonormal matrix and \(T_{1:i}\) is an upper triangular matrix.
        \IF{the stopping criterion is satisfied}
             \STATE Solve the triangular system \(T_{1:i} y^{(i)} = U_{1:i}\trans g\) to obtain \(y^{(i)}\).
             \RETURN \(x = x^{(i)} \gets x^{(0)} + Z_{1:i} y^{(i)}\).
        \ENDIF
    \ENDFOR
\end{algorithmic}
\end{algorithm}

\begin{algorithm}[H]
\begin{algorithmic}[1]
    \caption{Restarted preconditioned sketched GMRES \label{alg:restarted_sGMRES}}
    \REQUIRE
     A matrix \(A \in \mathbb R^{n\times n}\), a right-hand side \(b \in \mathbb R^{n}\), an initial approximation \(x^{(0)} \in \mathbb R^{n}\), the maximal number of iterations \(m\), a subspace embedding \(S \in \mathbb R^{s\times n}\) (\(s>m\)), left and right preconditioners \(M_L\) and \(M_R\), and the maximum number of restarts \(\texttt{nrestarts}\).
    \ENSURE
    A computed solution \(x \in \mathbb R^{n}\) approximating the solution of \(A x = b\).
    
    \FOR{j = 1:\texttt{nrestarts}}
        \STATE Call preconditioned sGMRES with input $M_{L}$, $M_{R}$, $A$, $b$, and \(x^{(0)}\) to obtain $x$.
        \STATE Update the initial approximate solution for the next iteration, $x^{(0)} \gets x$.
    \ENDFOR
\end{algorithmic}
\end{algorithm}

It is clear from experimental observations made in \cite{NT2024, guttel2024sketch} that the basis condition number should not grow too large throughout the iteration in order to ensure numerical stability of sketched GMRES.
The non-orthogonal basis can be constructed using the $t$ truncated Arnoldi process, as summarized in Algorithm~\ref{alg:trunc_Arnoldi}, where each new Krylov vector is orthogonalized against the previous $t$ vectors. Alternatively, the sketch-and-select method \cite{guttel2024sketch}
uses the sketched basis $S B_i$ to select $t$ vectors, among all previously computed vectors, with which to perform the partial orthogonalization, with the aim of generating a basis with smaller condition number than that generated by truncated Arnoldi.


\begin{algorithm}[!tb]
\begin{algorithmic}[1]
    \caption{The \(i\)-th iteration of the truncated Arnoldi process \label{alg:trunc_Arnoldi}}
    \REQUIRE
    A matrix \(A \in \mathbb R^{n\times n}\), a truncated parameter \(t\),
    the basis \(B_{1:i-1}\) generated by the first \(i-1\) truncated Arnoldi iterations.
    \ENSURE
    The basis \(B_{1:i}\).

    \STATE \(B_i = AB_{i-1}\).
    \FOR{\(l = \max(1, i-t+1):i\)} \label{line:algo_trunc:forloop}
        \STATE \(B_i = B_i - B_l(B_l\trans B_i)\).
    \ENDFOR
    \STATE \(B_i = B_i/\norm{B_i}\).
\end{algorithmic}
\end{algorithm}

\section{Backward stability of preconditioned sketched GMRES}
\label{sec:stability}
In this section, we aim to analyze the backward stability of preconditioned sketched GMRES (Algorithm~\ref{alg:sGMRES}) by primarily following the analysis of the MOD-GMRES framework~\cite{BHMV2024}.
The first \(i\) iterations of preconditioned sketched GMRES can be summarized in the following steps:

\begin{enumerate}
    \item Generate the preconditioned Krylov basis \(Z_{1:i} = M_R^{-1} B_{1:i}\) with the (non-orthogonal) Krylov basis \(B_{1:i}\), and compute \(W_{1:i} = M_L^{-1} A Z_{1:i}\). \label{eq:line-C}
    \item Solve the least squares problem 
    $ y^{(i)} = \argmin_{y} \norm{Sb - SW_{1:i} y}
    $ by solving a triangular linear system \eqref{eq:triangular_system} using the QR factorization of $C_{1:i}=SW_{1:i}$. \label{eq:line-ls}
    \item Compute the approximate solution \(x^{(i)} = Z_{1:i} y^{(i)}\). \label{eq:line-x}
\end{enumerate}
Notice that the Krylov basis \(B_{1:i}\) is non-orthogonal and may even be ill-conditioned, and the least squares problems solved in Step~\ref{eq:line-ls} are different for sketched GMRES and MOD-GMRES.

The main result of our backward error analysis relies on various assumptions on each of these steps, which we list here.
\paragraph{Step~\ref{eq:line-C}: Matrix--matrix product}
For the computed preconditioned basis \(\hat Z_{1:i} \in \mathbb{R}^{n\times i}\), we assume that the computed matrix--matrix product \(\hat W_{1:i}\) computed by Line~\ref{algo-line:computeW} in Algorithm~\ref{alg:sGMRES} satisfies
\begin{equation} \label{eq:compute-W}
    \hat W_{1:i} = M_L^{-1} A \hat Z_{1:i} + \Delta W_{1:i}, \quad \normF{\Delta W_{1:i}} \leq \epsW \norm{M_L^{-1}} \normF{A} \normF{\hat Z_{1:i}},
\end{equation}
with \(\epsW\leq \bigO(\macheps)\) from the standard rounding error analysis of~\cite{H2002}.
Similarly, we assume that the computed result of \(\tilde b := M_L^{-1}b\) computed by Line~\ref{algo-line:computer} in Algorithm~\ref{alg:sGMRES} satisfies
\begin{equation} \label{eq:compute-MLb}
    \hat b = \tilde b + \Delta b, \quad \norm{\Delta b}\leq \epsb \norm{M_L^{-1}} \norm{b}.
\end{equation}
Note that we assume \(x^{(0)} = 0\) in the analysis, and thus \(r^{(0)} = b\).

\paragraph{Step~\ref{eq:line-ls}: Least squares solver}
In iteratively solving the least squares problem in Step~\ref{eq:line-ls}, we assume that there is a special dimension \(k \leq n\), known as the \textit{key dimension} \cite{BHMV2024}, for which the computed solution is guaranteed to have reached a small backward error, i.e., \(\hat b\) lies in the range of non-singular basis \(\hat W_{1:k}\).
The definition of the key dimension is as follows:
\begin{align}
    \sigmin(\bmat{\phi \hat b& \hat W_{1:k}}) &\leq \epskapparW \normF{\bmat{\phi \hat b& \hat W_{1:k}}}, \label{eq:key-dim-rW} \\
    \sigmin(\hat W_{1:k}) &\geq \epskappaW \normF{\hat W_{1:k}} \label{eq:key-dim-W}
\end{align}
for all \(\phi > 0\), where \(\epskapparW\) and \(\epskappaW\) are scalar functions of \(n\), \(k\), and \(\macheps\), depending on the method for solving the least squares problem.
If \(\epskapparW=\bigO(\macheps)\) and \(\epskappaW=\bigO(\macheps)\), \eqref{eq:key-dim-W} guarantees that \(\hat W_{1:k}\) is a numerically non-singular basis of the subspace with dimension \(k\) and then \eqref{eq:key-dim-rW} means that \(\hat b\) lies in the range of \(\hat W_{1:k}\) numerically.
Furthermore, it also implies that \(\min_y \norm{\hat b - \hat W_{1:k} y}\) is sufficiently small.

Using the notion of the key dimension, we now aim to write down an expression for the computed sketched residual vector at the $k^{th}$ iteration of sketched GMRES. We then provide a bound on this sketched residual in Lemma \ref{lem:rkS-norm}, and use this bound in the proof of our final backward stability result in Theorem \ref{thm:stability_restarted_sgmres}.

We first set up the least squares problem by constructing the matrix $\hat{C}_{1:i}$, and the vector $g: = S M_L^{-1} b$. We assume that
\begin{align}
    \hat C_{1:i} &= S \hat W_{1:i} + \Delta C_{1:i}, \quad \norm{\Delta C_{j}} \leq \epsapplyS \norm{S \hat W_j}, \quad \forall j \leq i,  \label{eq:compute-C} \\
    \hat g &= S \hat b + \Delta g, \quad \norm{\Delta g} \leq \epsapplyS \norm{M_L^{-1}} \norm{S \hat b},  \label{eq:compute-g}
\end{align}
and furthermore, that the computed solution \(\hat y^{(i)}\) of the least squares problem in Step~\ref{eq:line-ls} satisfies
\begin{equation} \label{eq:sls}
    \begin{split}
        &\hat y^{(i)} = \argmin_y \norm{\hat g + \Delta \hat g - (\hat C_{1:i} + \Delta \hat C_{1:i}) y}, \\
        &\norm{\bmat{\Delta \hat g& \Delta \hat C_{1:i}} e_j} \leq \epssls \norm{\bmat{\hat g& \hat C_{1:i}} e_j}, \quad \forall j \leq i + 1,
    \end{split}
\end{equation}
with \(\epssls \leq \bigO(\macheps)\) from~\cite[Theorem 20.3]{H2002}.
Substituting \(\hat C_{1:i}\) and \(\hat g\) in~\eqref{eq:sls} with~\eqref{eq:compute-C} and~\eqref{eq:compute-g}, we have
\begin{equation} \label{eq:ls-S}
    \begin{split}
        & \hat y^{(i)} = \argmin_y \normbig{S\tilde \hat b + \Delta \tilde g - (S \hat W_{1:i} + \Delta \tilde C_{1:i}) y}, \\
        & \norm{\Delta \tilde g} \leq \norm{\Delta g} + \norm{\Delta\hat g} \leq (1 + \epssls) \epsapplyS \norm{S\hat b} + \epssls \norm{S\hat b} \leq \bigl(\epssls + (1 + \epssls) \epsapplyS\bigr) \norm{S\hat b}, \\
        & \norm{\Delta \tilde C_{j}} \leq \norm{\Delta C_j} + \norm{\Delta\hat C_j} \leq \bigl(\epssls + (1 + \epssls) \epsapplyS\bigr) \norm{S \hat W_{j}}, \quad \forall j \leq i.
    \end{split}
\end{equation}

We can now prove a bound on the norm of the sketched residual at the key dimension $k$, i.e., the norm of the vector
\begin{equation} \label{eq:def-rkS}
r^{(k)}_S := S\hat b + \Delta \tilde g - (S \hat W_{1:k} + \Delta \tilde C_{1:k}) \hat y^{(k)}.
\end{equation}

\begin{lemma} \label{lem:rkS-norm}
    Assume that \(\hat W_{1:k}\) satisfies
    \begin{align}
        \sigmin(\bmat{\phi S \hat b& S \hat W_{1:k}}) &\leq c_1(n, k) \bigl(\epsapplyS + \epssls\bigr) \normF{\bmat{\phi S \hat b& S \hat W_{1:k}}}, \label{eq:lem:assump:SbW} \\
        \sigmin(S \hat W_{1:k}) &\geq c_2(n, k)\bigl(\epsapplyS + \epssls\bigr) \normF{S \hat W_{1:k}}, \label{eq:lem:assump:SW}
    \end{align}
    where $S$ is an $\epsilon$-subspace embedding for the space spanned by the matrix \(\bmat{\hat b& \hat W_{1:k}}\) satisfying~\eqref{eq:sketch}, and \(c_1(n, k)\) and \(c_2(n, k)\)
    are low degree polynomials related to \(n\) and \(k\).
    Also assume that \(\hat y^{(k)}\) satisfies~\eqref{eq:ls-S}.
    Then
    \begin{equation} \label{eq:rkS-norm}
        \norm{r^{(k)}_S} \leq c(n, k) \bigl(\epsapplyS + \epssls\bigr) \normF{S \hat W_{1:k}} \norm{\hat y^{(k)}},
    \end{equation}
    where \(c(n, k)\) is also a low degree polynomial related to \(n\) and \(k\).
\end{lemma}

\begin{proof}
    We follow the proof in~\cite{BHMV2024} and  employ~\cite[Theorem~2.4]{paige2006modified} to bound \(\norm{r^{(k)}_S}\).
    Applying~\cite[Theorem~2.4]{paige2006modified} to the least squares problem~\eqref{eq:ls-S}, we have
    \begin{equation} \label{eq:lem-rkS-norm-proof:norm-rkS-0}
        \norm{r^{(k)}_S}^2\leq \sigmin^2\bigl(\bmat{\phi(S\hat b + \Delta \tilde g) & S \hat W_{1:k} + \Delta \tilde C_{1:k}}\bigr) \biggl(\frac{1}{\phi^{2}}+\frac{\norm{\hat y^{(k)}}^2}{1-\delta_k^2}\biggr),
    \end{equation}
    with \(\delta_k\) defined by
    \begin{equation} \label{eq:lem-rkS-norm-proof:defdelta}
        \delta_k = \frac{\sigmin\bigl(\bmat{\phi(S\hat b + \Delta \tilde g) & S \hat W_{1:k} + \Delta \tilde C_{1:k}}\bigr)}{\sigmin\bigl( S \hat W_{1:k} + \Delta \tilde C_{1:k}\bigr)} < 1.
    \end{equation}
    We choose \(\phi\) such that \(\phi^{-2} = \norm{\hat y^{(k)}}^2/(1-\delta_k^2)\) as in~\cite[(3.17)]{BHMV2024}, which can simplify~\eqref{eq:lem-rkS-norm-proof:norm-rkS-0} as
    \begin{equation} \label{eq:lem-rkS-norm-proof:norm-rkS-1}
        \norm{r^{(k)}_S}\leq \sqrt{2}\,\sigmin\bigl(\bmat{\phi(S\hat b + \Delta \tilde g) & S \hat W_{1:k} + \Delta \tilde C_{1:k}}\bigr) \phi^{-1}.
    \end{equation}
    Then we will show that there exists a \(\phi>0\) satisfying \(\phi^{-2} = \norm{\hat y^{(k)}}^2/(1-\delta_k^2)\).
    This is equivalent to showing that there exists a \(\phi>0\) such that \(f(\phi)=0\), where the function \(f(\phi)\) is defined as
    \begin{equation}
    \begin{split}
        f(\phi) &= \sigmin^2\bigl( S \hat W_{1:k} + \Delta \tilde C_{1:k}\bigr) - \sigmin^2\bigl(\bmat{\phi(S\hat b + \Delta \tilde g) & S \hat W_{1:k} + \Delta \tilde C_{1:k}}\bigr) \\
        &\quad- \norm{\hat y^{(k)}}^2\phi^{2} \sigmin^2\bigl( S \hat W_{1:k} + \Delta \tilde C_{1:k}\bigr).
    \end{split}
    \end{equation}
    Notice that \(f(0)> 0\) and \(f(\norm{\hat y^{(k)}}^{-1})< 0\).
    By continuity, this implies that there exists \(\phi\in (0, \norm{\hat y^{(k)}}^{-1})\) such that \(\phi^{-2} = \norm{\hat y^{(k)}}^2/(1-\delta_k^2)\) and \(\delta_k<1\).

    It remains to bound \(\sigmin\bigl(\bmat{\phi(S\hat b + \Delta \tilde g) & S \hat W_{1:k} + \Delta \tilde C_{1:k}}\bigr)\) and \(\phi\) in~\eqref{eq:lem-rkS-norm-proof:norm-rkS-1}.
    By the assumption~\eqref{eq:lem:assump:SbW}, and the estimation of \(\norm{\Delta \tilde g}\) and \(\norm{\Delta \tilde C_{1:k}}\) shown in~\eqref{eq:ls-S}, we have
    \begin{equation} \label{eq:lem-rkS-norm-proof:sigmin}
    \begin{split}
        \sigmin\bigl(&\bmat{\phi(S\hat b + \Delta \tilde g) & S \hat W_{1:k} + \Delta \tilde C_{1:k}}\bigr) \\
        &\leq c_1(n, k) \bigl(\epsapplyS + \epssls\bigr) \normF{\bmat{\phi S \hat b& S \hat W_{1:k}}} \\
        &\leq c_1(n, k) \bigl(\epsapplyS + \epssls\bigr) \bigl(\norm{S \hat b}\phi + \normF{S \hat W_{1:k}}\bigr).
    \end{split}
    \end{equation}
    To establish a bound for \(\phi\), it suffices to bound \(\delta_k\).
    Combining~\eqref{eq:lem-rkS-norm-proof:sigmin} with the assumption~\eqref{eq:lem:assump:SW} and the estimation~\eqref{eq:ls-S} of \(\norm{\Delta \tilde C_{1:k}}\), \(\delta_k\) can be bounded by
    \begin{equation} \label{eq:lem-rkS-norm-proof:deltakbound}
        \begin{split}
            \delta_k \leq \frac{c_1(n, k) \bigl(\epsapplyS + \epssls\bigr) \bigl(\norm{S \hat b}\phi + \normF{S \hat W_{1:k}}\bigr)}{c_2(n, k) \bigl(\epsapplyS + \epssls\bigr) \normF{S \hat W_{1:k}}}.
        \end{split}
    \end{equation}
    We now show that \(\norm{S \hat b}\phi\) can be bounded by \(\normF{S \hat W_{1:k}}\).
    From the definition of \(r^{(k)}_S\) in~\eqref{eq:def-rkS}, we have
    \begin{equation*}
        \norm{S \hat b}\phi \leq \norm{r^{(k)}_S}\phi + \norm{\Delta \tilde g}\phi
        + \normF{S \hat W_{1:k}} \norm{\hat y^{(k)}}\phi + \norm{\Delta\tilde C_{1:k}} \norm{\hat y^{(k)}}\phi.
    \end{equation*}
    Notice that \(\norm{\hat y^{(k)}}\phi\leq 1\) since \(\phi\in (0, \norm{\hat y^{(k)}}^{-1})\).
    Furthermore, using~\eqref{eq:lem-rkS-norm-proof:norm-rkS-1}, \eqref{eq:lem-rkS-norm-proof:sigmin}, and the estimation of \( \norm{\Delta \tilde g}\) and \(\norm{\Delta \tilde C_{1:k}}\) from~\eqref{eq:ls-S}, we obtain
    \begin{equation*}
        \norm{S \hat b}\phi \leq \tilde c_1(n, k)(\epssls+\epsapplyS) \bigl(\norm{S \hat b}\phi + \normF{S \hat W_{1:k}}\bigr) + \normF{S \hat W_{1:k}},
    \end{equation*}
    which means
    \begin{equation} \label{eq:lem-rkS-norm-proof:norm-Sb}
        \norm{S \hat b}\phi \leq \frac{1+\tilde c_1(n, k)(\epssls+\epsapplyS)}{1-\tilde c_1(n, k)(\epssls+\epsapplyS)} \normF{S \hat W_{1:k}},
    \end{equation}
    and then \(\delta_k\leq 1/2\) by choosing appropriate \(c_1(n, k)\) and \(c_2(n, k)\) in~\eqref{eq:lem-rkS-norm-proof:deltakbound}.
    Thus, we have
    \begin{equation} \label{eq:lem-rkS-norm-proof:phibound}
        \phi = \frac{\norm{\hat y^{(k)}}}{\sqrt{1-\delta_k}}
        \leq 2\norm{\hat y^{(k)}},
    \end{equation}
    and, using \eqref{eq:lem-rkS-norm-proof:norm-Sb}, \eqref{eq:lem-rkS-norm-proof:phibound}, \eqref{eq:lem-rkS-norm-proof:sigmin}, and~\eqref{eq:lem-rkS-norm-proof:norm-rkS-1}, we can conclude the proof.
\end{proof}

Lemma~\ref{lem:rkS-norm} provides the values of \(\epskapparW\) and \(\epskappaW\) for defining the key dimension, that is, \(c_1(n, k) \bigl(\epsapplyS + \epssls\bigr)\) and \(c_2(n, k) \bigl(\epsapplyS + \epssls\bigr)\), respectively.

\paragraph{Step~\ref{eq:line-x}: Computation of the approximate solution}
For \(\hat Z_{1:i} \in \mathbb{R}^{n\times i}\) and \(\hat y^{(i)} \in \mathbb{R}^i\), we assume that
\begin{equation} \label{eq:updatex}
    \hat x^{(i)} = \hat Z_{1:i} \hat y^{(i)} + \Delta x^{(i)}, \quad \norm{\Delta x^{(i)}} \leq \epsx \norm{\hat Z_{1:i}} \norm{\hat y^{(i)}}
\end{equation}
with \(\epsx \leq \bigO(\macheps)\) from~\cite{H2002}.

\medskip
We are now ready to describe the backward stability of preconditioned sketched GMRES (Algorithm~\ref{alg:sGMRES}) in the following theorem.

\begin{theorem}\label{thm:stability_sgmres}
    Assume that Algorithm~\ref{alg:sGMRES} is applied with a computed preconditioned basis \(\hat Z_{1:m} \in \mathbb{R}^{n\times m}\) and a randomized subspace embedding \(S \in \mathbb{R}^{s\times n}\) satisfying~\eqref{eq:sketch} for the subspace \(\mathbb{R}^{n\times m}\) with high probability, where~\eqref{eq:compute-W}, \eqref{eq:compute-MLb}, \eqref{eq:compute-C}, \eqref{eq:compute-g}, \eqref{eq:sls}, and~\eqref{eq:updatex} are satisfied.
    If there exists \(k \leq m\) such that \(c_1(n, k) \bigl(\epsapplyS + \epssls\bigr) \sqrt{\frac{1 - \epsilon}{1 + \epsilon}} \kappa(\bmat{\phi \hat b& \hat W_{1:k}}) \geq 1\),
    \begin{equation} \label{eq:thm-assump}
        \sqrt{k}\, c_2(n, k)\bigl(\epsapplyS + \epssls\bigr) \sqrt{\frac{1 + \epsilon}{1 - \epsilon}} \kappa(\hat W_{1:k}) \leq \frac{1}{2}, \quad\text{and}\quad
        \epsx \kappa(\hat Z_{1:k}) \leq \frac{1}{2},
    \end{equation}
    then, with high probability,
    \begin{equation} \label{eq:res-xk}
        \begin{split}
            \lVert b - &A \hat x^{(k)}\rVert \\
            &\leq \tilde c(n, k) \biggl(\frac{1 + \epsilon}{1 - \epsilon} \bigl(\epsapplyS + \epssls\bigr) + \epsW + \epsx\biggr) \kappa(M_L) \kappa(\hat Z_{1:k}) \normF{A} \norm{\hat x^{(k)}} \\
            &\quad + \biggl(\frac{1 + \epsilon}{1 - \epsilon} (1+\epsb)\bigl(\epssls + (1 + \epssls) \epsapplyS\bigr) + \epsb \biggr) \kappa(M_L) \norm{b}
        \end{split}
    \end{equation}
    with a low degree polynomial \(\tilde c(n, k)\) related to \(n\) and \(k\).
    Furthermore, if \(\epsW \), \(\epsb\), \(\epssls\), \(\epsapplyS\), \(\epsx \leq \bigO(\macheps)\), then
    \begin{equation} \label{eq:thm:res-xk-simple}
        \norm{b - A \hat x^{(k)}} \leq \frac{1 + \epsilon}{1 - \epsilon} \bigO(\macheps)\kappa(M_L) \bigl(\kappa(\hat Z_{1:k}) \normF{A} \norm{\hat x^{(k)}} + \norm{b}\bigr).
    \end{equation}
\end{theorem}

\begin{proof}
    The residual of \(\hat x^{(k)}\) can be bounded by the preconditioned residual as follows:
    \begin{equation} \label{eq:thm-proof:b-Axk-0}
        \norm{b - A \hat x^{(k)}} = \norm{M_L\bigl(M_L^{-1} b - M_L^{-1} A \hat x^{(k)})\bigr}
        \leq \norm{M_L} \norm{M_L^{-1} b - M_L^{-1} A \hat x^{(k)}},
    \end{equation}
    which means that we only need to focus on the preconditioned residual \(\norm{M_L^{-1} b - M_L^{-1} A \hat x^{(k)}}\).
    By substituting \(\hat x^{(k)}\) with~\eqref{eq:updatex}, \(M_L^{-1}b\) with~\eqref{eq:compute-MLb}, and \(M_L^{-1} A \hat Z_{1:k}\) with~\eqref{eq:compute-W}, we have
    \begin{equation} \label{eq:thm-proof:b-Axk}
        \norm{M_L^{-1} b - M_L^{-1} A \hat x^{(k)}}
        \leq \norm{\hat b - \hat W_{1:k} \hat y^{(k)}} + \norm{\Delta W_{1:k} \hat y^{(k)}} + \norm{\Delta b} + \norm{M_L^{-1}} \norm{A \Delta x^{(k)}}.
    \end{equation}
    Furthermore, from~\cite[Fact 2.2]{NT2024} and the definition~\eqref{eq:def-rkS} of \(r^{(k)}_S\), we have
    \begin{equation} \label{eq:thm-proof:b-Wy}
    \begin{split}
        \norm{\hat b - \hat W_{1:k} \hat y^{(k)}}
        &\leq \frac{1}{1 - \epsilon} \norm{S \hat b - S \hat W_{1:k} \hat y^{(k)}} \\
        &\leq \frac{1}{1 - \epsilon} \norm{r^{(k)}_S} + \frac{1}{1 - \epsilon} \bigl(\norm{\Delta \tilde g} + \normF{\Delta \tilde C_{1:k}} \norm{\hat y^{(k)}}\bigr).
    \end{split}
    \end{equation}
    Combining~\eqref{eq:thm-proof:b-Wy} with~\eqref{eq:thm-proof:b-Axk}, it follows that
    \begin{equation} \label{eq:thm:proof:res}
        \begin{split}
            \norm{M_L^{-1} \hat b - M_L^{-1} A \hat x^{(k)}}
            &\leq \frac{1}{1 - \epsilon} \norm{r^{(k)}_S} + \frac{1}{1 - \epsilon} \bigl(\norm{\Delta \tilde g} + \normF{\Delta \tilde C_{1:k}} \norm{\hat y^{(k)}} \bigr) \\
            &\quad + \norm{\Delta W_{1:k}} \norm{\hat y^{(k)}} + \norm{\Delta b}
            + \norm{M_L^{-1}} \normF{A} \norm{\Delta x^{(k)}}.
        \end{split}
    \end{equation}

    It remains to bound \(\norm{r^{(k)}_S}\) and \(\norm{\hat y^{(k)}}\).
    To apply Lemma~\ref{lem:rkS-norm} to bound \(\norm{r^{(k)}_S}\), we first need to verify~\eqref{eq:lem:assump:SbW} and~\eqref{eq:lem:assump:SW}.
    From~\cite[Corollary 2.2]{BG2022}, the condition numbers of the sketched matrices satisfy
    \begin{equation*}
        \kappa(S \bmat{\phi \hat b& \hat W_{1:k}}) \geq \sqrt{\frac{1 - \epsilon}{1 + \epsilon}} \kappa(\bmat{\phi \hat b& \hat W_{1:k}}) \quad \text{and} \quad
        \kappa(S \hat W_{1:k}) \leq \sqrt{\frac{1 + \epsilon}{1 - \epsilon}} \kappa(\hat W_{1:k}),
    \end{equation*}
    which implies, from \(c_1(n, k) \bigl(\epsapplyS + \epssls\bigr) \sqrt{\frac{1 - \epsilon}{1 + \epsilon}} \kappa(\bmat{\phi b& \hat W_{1:k}}) \geq 1\) and~\eqref{eq:thm-assump},
    \begin{align*}
        \sigmin(S \bmat{\phi \hat b& \hat W_{1:k}}) &\leq c_1(n, k) \bigl( \epsapplyS + \epssls\bigr) \norm{S \bmat{\phi \hat b& \hat W_{1:k}}} \\
        &\leq c_1(n, k) \bigl( \epsapplyS + \epssls\bigr) \normF{S \bmat{\phi b& \hat W_{1:k}}}, \\
        \sqrt{k}\, c_2(n, k) \bigl(\epsapplyS + \epssls\bigr) \kappa(S \hat W_{1:k}) &\leq \sqrt{k}\, c_2(n, k) \bigl(\epsapplyS + \epssls\bigr) \sqrt{\frac{1 + \epsilon}{1 - \epsilon}} \kappa(\hat W_{1:k}) \leq 1,
    \end{align*}
    which guarantees that~\eqref{eq:lem:assump:SbW} and~\eqref{eq:lem:assump:SW} hold.
    Thus, from Lemma~\ref{lem:rkS-norm}, we obtain the bound on \(\norm{r^{(k)}_S}\):
    \begin{equation} \label{eq:thm-proof:rkSnorm}
        \norm{r^{(k)}_S} \leq c(n, k) \bigl(\epsapplyS + \epssls\bigr) (1+\epsilon) \normF{\hat W_{1:k}} \norm{\hat y^{(k)}}.
    \end{equation}
    Using the bounds on \(\norm{r^{(k)}_S}\), i.e.,~\eqref{eq:thm-proof:rkSnorm}, as well as the bounds on \(\Delta \tilde g\), \(\Delta \tilde C_{1:k}\), \(\Delta W_{1:k}\), \(\Delta b\), and \(\Delta x^{(k)}\) shown in~\eqref{eq:ls-S}, \eqref{eq:compute-W}, \eqref{eq:compute-MLb}, and~\eqref{eq:updatex}, then \eqref{eq:thm:proof:res} can be bounded as
    \begin{equation} \label{eq:thm:proof:res-2}
        \begin{split}
            \norm{M_L^{-1} \hat b - M_L^{-1} A \hat x^{(k)}}
            &\leq \frac{1 + \epsilon}{1 - \epsilon} c(n, k) \bigl( \epsapplyS + \epssls\bigr) \normF{\hat W_{1:k}} \norm{\hat y^{(k)}} \\
            &\quad + \frac{1 + \epsilon}{1 - \epsilon} \bigl(\epssls + (1 + \epssls) \epsapplyS\bigr) \biggl(\norm{\hat b} + \normF{\hat W_{1:k}} \norm{\hat y^{(k)}} \biggr) \\
            &\quad + \epsW \norm{M_L^{-1}} \normF{A} \normF{\hat Z_{1:k}} \norm{\hat y^{(k)}} + \epsb \norm{M_L^{-1}} \norm{b} \\
            &\quad + \epsx \norm{M_L^{-1}}  \normF{A} \normF{\hat Z_{1:k}}  \norm{\hat y^{(k)}}.
        \end{split}
    \end{equation}
    Furthermore, noticing that \(\normF{\hat W_{1:k}} \leq (1+\epsW) \norm{M_L^{-1}} \normF{A} \normF{\hat Z_{1:k}}\) from~\eqref{eq:compute-W}, \eqref{eq:thm:proof:res-2} can be simplified as
    \begin{equation} \label{eq:thm:proof:res-3}
        \begin{split}
            \lVert &M_L^{-1} \hat b - M_L^{-1} A \hat x^{(k)}\rVert \\
            &\leq \biggl(\frac{1 + \epsilon}{1 - \epsilon} c(n, k) (1+\epsW) \bigl(\epsapplyS + \epssls\bigr) + \sqrt{k}\,\epsW + \sqrt{k}\,\epsx\biggr) \norm{M_L^{-1}} \normF{A} \norm{\hat Z_{1:k}} \norm{\hat y^{(k)}} \\
            &\quad + \biggl(\frac{1 + \epsilon}{1 - \epsilon} (1+\epsb)\bigl(\epssls + (1 + \epssls) \epsapplyS\bigr) + \epsb \biggr) \norm{M_L^{-1}} \norm{b}.
        \end{split}
    \end{equation}

    Thus, using~\eqref{eq:updatex} we have
    \begin{equation}\label{eq:ybound}
        \norm{\hat y^{(k)}} \leq \frac{\norm{\hat x^{(k)}}}{\sigmin(\hat Z_{1:k})- \epsx \norm{\hat Z_{1:k}}}.
    \end{equation}
    Combining~\eqref{eq:thm:proof:res-3} with~\eqref{eq:ybound}, we derive
    \begin{equation*} \label{eq:thm:proof:res-4}
        \begin{split}
            \lVert &M_L^{-1} \hat b - M_L^{-1} A \hat x^{(k)}\rVert \\
            &\leq \biggl(\frac{1 + \epsilon}{1 - \epsilon} c(n, k)(1+\epsW)  \bigl(\epsapplyS + \epssls\bigr) + \sqrt{k}\,\epsW + \sqrt{k}\,\epsx\biggr) \frac{ \norm{M_L^{-1}} \normF{A} \norm{\hat Z_{1:k}} \norm{\hat x^{(k)}}}{\sigmin(\hat Z_{1:k})- \epsx \norm{\hat Z_{1:k}}} \\
            &\quad + \biggl(\frac{1 + \epsilon}{1 - \epsilon} (1+\epsb)\bigl(\epssls + (1 + \epssls) \epsapplyS\bigr) + \epsb \biggr) \norm{M_L^{-1}} \norm{b} \\
            &\leq \biggl(\frac{1 + \epsilon}{1 - \epsilon} c(n, k)(1+\epsW)  \bigl(\epsapplyS + \epssls\bigr) + \sqrt{k}\,\epsW + \sqrt{k}\,\epsx\biggr) \frac{\norm{M_L^{-1}} \normF{A} \kappa(\hat Z_{1:k}) \norm{\hat x^{(k)}}}{1 - \epsx \kappa(\hat Z_{1:k})} \\
            &\quad + \biggl(\frac{1 + \epsilon}{1 - \epsilon} (1+\epsb)\bigl(\epssls + (1 + \epssls) \epsapplyS\bigr) + \epsb \biggr) \norm{M_L^{-1}} \norm{b},
        \end{split}
    \end{equation*}
    which, together with~\eqref{eq:thm-proof:b-Axk-0} and the assumption~\eqref{eq:thm-assump} concludes the proof.
\end{proof}

\begin{remark}
    Theorem~\ref{thm:stability_sgmres} indicates that~\eqref{eq:res-xk} and~\eqref{eq:thm:res-xk-simple} are satisfied with high probability.
    This is due to the choice of the subspace embedding \(S\), which generally meets the condition~\eqref{eq:sketch} with high probability.
    Practically, the sketching method rarely fails, and slightly increasing the embedding dimension can further improve its reliability.
    In the worst case, \(\epsilon\) may approach \(1\) and the constant \((1+\epsilon)/(1-\epsilon)\) may be relatively large.
    For this reason, it is distracting to quantify the failure probabilities; see details in~\cite{MT2020}.
\end{remark}

\begin{remark}
    For standard GMRES, the existence of such a \(k\) is straightforward to be guaranteed by observing that \(\hat b \in \mathbb{R}^n\), because the condition
    \(
    c_1(n, k)\bigl(\epsapplyS + \epssls\bigr)\sqrt{\frac{1 - \epsilon}{1 + \epsilon}}\, \kappa(\bmat{\phi \hat b & \hat W_{1:k}}) \geq 1
    \)
    is automatically satisfied when \(k = n\).
    In contrast, for sketched GMRES, one must select a proper subspace embedding \(S \in \mathbb{R}^{s \times n}\) for the \(m\)-dimensional subspace in order to ensure the existence of such a \(k\).
    When \(s\) is sufficiently large—particularly in the extreme case \(s = n\)—the existence of \(k\) can be guaranteed in essentially the same way as for standard GMRES, but then the advantage of employing random sketching is lost.
\end{remark}

\begin{remark} \label{remark:thm1-res}
    Let \(\varepsilon_*\) denote all the \(\varepsilon\) with any subscript mentioned in this section.
    Assuming \(\varepsilon_*\leq \bigO(\macheps)\), we simplify Theorem~\ref{thm:stability_sgmres} for the case without preconditioning, i.e., \(M_L = M_R = I\) and \(\kappa(\hat{Z}_{1:k}) = \kappa(\hat{B}_{1:k})\): if it holds for the key dimension \(k\) that
    \begin{equation}
        \bigO(\macheps) \sqrt{\frac{1 + \epsilon}{1 - \epsilon}} \kappa(\hat{W}_{1:k})\leq 1/2, \qquad
        \bigO(\macheps) \kappa(\hat{B}_{1:k}) \leq 1/2,
    \end{equation}
    then the relative backward error can be bounded as
    \begin{equation}\label{eq:remark1_bound} 
        \frac{\norm{b - A \hat x^{(k)}}}{\normF{A} \norm{\hat x^{(k)}} + \norm{b}} \leq \frac{1 + \epsilon}{1 - \epsilon} \bigO(\macheps) \kappa(\hat B_{1:k}).
    \end{equation}
\end{remark}

\begin{remark}
    Considering \(\varepsilon_*\leq \bigO(\macheps)\), note that the assumption~\eqref{eq:thm-assump} in Theorem \ref{thm:stability_sgmres} can be replaced by
    \begin{equation*}\bigO(\macheps) \sqrt{\frac{1 + \epsilon}{1 - \epsilon}} \kappa(M_L) \kappa(A) \kappa(\hat B_{1:k}) \leq 1/2\end{equation*} 
    since
    \begin{equation} \label{eq:remark:kappaW}
        \begin{split}
            \kappa(\hat W_{1:k}) &\leq \frac{\norm{M_L^{-1} A \hat Z_{1:k}} + \bigO(\macheps) \norm{M_L^{-1}} \normF{A} \norm{\hat B_{1:k}}}{\sigmin(M_L^{-1} A \hat Z_{1:k}) - \bigO(\macheps) \norm{M_L^{-1}} \normF{A} \norm{\hat Z_{1:k}}} \\
            &\leq 2\kappa(M_L^{-1} A \hat Z_{1:k}) + \bigO(\macheps)\kappa(M_L) \kappa(A) \kappa(\hat Z_{1:k}) \\
            &\leq (2 + \bigO(\macheps))\kappa(M_L) \kappa(A) \kappa(\hat Z_{1:k})
        \end{split}
    \end{equation}
    from~\eqref{eq:compute-W}.
    Therefore, if \(S\) is chosen to be a \(n\times n\) identity matrix and \(\hat B_{1:k}\) is chosen to be an orthonormal Krylov basis, Theorem~\ref{thm:stability_sgmres} recovers a similar backward stability result as~\cite[Theorem~3.1]{BHMV2024}.
\end{remark}

Notice that \(\kappa(\hat W_{1:k})\leq 3 \kappa(M_L^{-1} A Z_{1:j})\) from~\eqref{eq:remark:kappaW} if \(\bigO(\macheps)\kappa(M_L) \kappa(A) \kappa(\hat Z_{1:k})\leq 1\).
As a consequence of Theorem \ref{thm:stability_sgmres}, we can conclude that if $\kappa(M_L^{-1} A Z_{1:j})$ is sufficiently small, then the backward stability of sGMRES depends on the condition number of the preconditioned basis $Z_{1:j}$. In particular, sGMRES is backward stable provided $\kappa(Z_{1:j})$ is not too large. 

It has been observed in numerical experiments in~\cite[Figure~9]{guttel2024sketch} that in some cases, sGMRES will continue to converge even when \(\kappa(\hat Z_{1:k})>10^{15}\), and so the bound shown in Theorem~\ref{thm:stability_sgmres} is not tight.
Therefore, we present the following lemma to show a tighter bound, which can be directly derived from~\eqref{eq:thm:proof:res-3} by noticing
\begin{equation*}
    \norm{b - A \hat x^{(k)}} \leq \norm{M_L} \lVert M_L^{-1} \hat b - M_L^{-1} A \hat x^{(k)}\rVert \quad\text{and}\quad
    \norm{M_L} \norm{M_L^{-1}} = \kappa(M_L).
\end{equation*}

\begin{lemma} \label{lem:backward-sGMRES}
    Assume that Algorithm~\ref{alg:sGMRES} is applied with a computed preconditioned basis \(\hat Z_{1:m} \in \mathbb{R}^{n\times m}\) and a randomized subspace embedding \(S \in \mathbb{R}^{s\times n}\) satisfying~\eqref{eq:sketch} for the subspace \(\mathbb{R}^{n\times m}\) with high probability, where~\eqref{eq:compute-W}, \eqref{eq:compute-MLb}, \eqref{eq:compute-C}, \eqref{eq:compute-g}, \eqref{eq:sls}, and~\eqref{eq:updatex} are satisfied.
    If there exists \(k \leq m\) such that \(c_1(n, k) \bigl(\epsapplyS + \epssls\bigr) \sqrt{\frac{1 - \epsilon}{1 + \epsilon}} \kappa(\bmat{\phi \hat b& \hat W_{1:k}}) \geq 1\),
    \begin{equation*}
        \sqrt{k}\, c_2(n, k)\bigl(\epsapplyS + \epssls\bigr) \sqrt{\frac{1 + \epsilon}{1 - \epsilon}} \kappa(\hat W_{1:k}) \leq \frac{1}{2}, \quad\text{and}\quad
        \epsx \kappa(\hat Z_{1:k}) \leq \frac{1}{2},
    \end{equation*}
    then, with high probability,
    \begin{equation} \label{eq:lem:res-xk}
        \begin{split}
            \norm{b - A \hat x^{(k)}}
            &\leq \tilde c(n, k) \biggl(\frac{1 + \epsilon}{1 - \epsilon} \bigl(\epsapplyS + \epssls\bigr) + \epsW + \epsx\biggr) \kappa(M_L) \norm{\hat Z_{1:k}} \normF{A} \norm{\hat y^{(k)}} \\
            &\quad + \biggl(\frac{1 + \epsilon}{1 - \epsilon} (1+\epsb)\bigl(\epssls + (1 + \epssls) \epsapplyS\bigr) + \epsb \biggr) \kappa(M_L) \norm{b}
        \end{split}
    \end{equation}
    with a low degree polynomial \(\tilde c(n, k)\) related to \(n\) and \(k\).
    Furthermore, if \(\epsW \), \(\epsb\), \(\epssls\), \(\epsapplyS\), \(\epsx \leq \bigO(\macheps)\), then
    \begin{equation} \label{eq:lem:res-xk-simple}
        \norm{b - A \hat x^{(k)}} \leq \frac{1 + \epsilon}{1 - \epsilon} \bigO(\macheps)\kappa(M_L) \bigl(\norm{\hat Z_{1:k}} \normF{A} \norm{\hat y^{(k)}} + \norm{b}\bigr).
    \end{equation}
\end{lemma}

\begin{remark}
Comparing the bounds given in~\eqref{eq:res-xk} and~\eqref{eq:thm:res-xk-simple} within Theorem~\ref{thm:stability_sgmres}, Lemma~\ref{lem:backward-sGMRES} utilizes \(\norm{\hat y^{(k)}}\) in the bounds~\eqref{eq:lem:res-xk} and~\eqref{eq:lem:res-xk-simple} as an alternative to employing \(\norm{\hat x^{(k)}}/\sigmin(\hat Z_{1:k})\).
Based on~\eqref{eq:ybound} and the presumption~\eqref{eq:thm-assump}, it follows that 
\[
\norm{\hat y^{(k)}}\leq \frac{2\norm{\hat x^{(k)}}}{\sigmin(\hat Z_{1:k})}.
\]
This implies that Lemma~\ref{lem:backward-sGMRES} offers a sharper bound compared to Theorem~\ref{thm:stability_sgmres}, as demonstrated by
\begin{equation*}
    \kappa(M_L) \bigl(\norm{\hat Z_{1:k}} \normF{A} \norm{\hat y^{(k)}} + \norm{b}\bigr)\leq 2\kappa(M_L) \bigl(\kappa(\hat Z_{1:k}) \normF{A} \norm{\hat x^{(k)}} + \norm{b}\bigr).
\end{equation*}
Furthermore, in practice, \(\norm{\hat y^{(k)}}\) can be much smaller than \(\norm{\hat x^{(k)}}/\sigmin(\hat Z_{1:k})\).
This implies that the relative backward error bound based on \(\norm{\hat y^{(k)}}\) can remain small or potentially decrease even when \(\kappa(\hat Z_{1:k})>10^{15} \).
\end{remark}

\begin{remark}
By Lemma~\ref{lem:backward-sGMRES}, the relative backward error can be bounded as
\begin{equation*}
    \frac{\norm{b - A \hat x^{(k)}}}{\normF{A} \norm{\hat x^{(k)}} + \norm{b}} \leq \frac{1 + \epsilon}{1 - \epsilon} \bigO(\macheps) \kappa(M_L) \frac{\norm{\hat Z_{1:k}} \norm{\hat y^{(k)}}}{\norm{\hat x^{(k)}}},
\end{equation*}
and can be simplified under the condition \(M_L=M_R=I\) to
\begin{equation} \label{eq:lem-res}
    \frac{\norm{b - A \hat x^{(k)}}}{\normF{A} \norm{\hat x^{(k)}} + \norm{b}} \leq \frac{1 + \epsilon}{1 - \epsilon} \bigO(\macheps) \frac{\norm{\hat B_{1:k}} \norm{\hat y^{(k)}}}{\norm{\hat x^{(k)}}}.
\end{equation}
\end{remark}

\section{Restarted sGMRES with adaptive truncation size}
\label{sec:restart-sGMRES}
Similar to the case of standard GMRES, restarting can be regarded as an iterative refinement process as follows:
\begin{enumerate}
    \item Compute the residual \(r^{(i)} = b - Ax^{(i)}\).
    \item Solve \(Ad^{(i)} = r^{(i)}\), which is computed by Algorithm~\ref{alg:sGMRES}.
    \item Update the solution \(x^{(i+1)} = x^{(i)} + d^{(i)}\).
\end{enumerate}
Furthermore, it can remove the influence of both \(\kappa(M_L)\) and $\kappa(Z_{1:j})$ on the backward stability of sGMRES, as shown in the following version of \cite[Theorem 4.1]{BHMV2024} describing the stability of restarted sketched GMRES (Algorithm \ref{alg:restarted_sGMRES}).

\begin{theorem}\label{thm:stability_restarted_sgmres}
    Assume that Algorithm~\ref{alg:sGMRES} is restarted every \(m\) iterations with a randomized subspace embedding \(S \in \mathbb{R}^{s\times n}\) satisfying~\eqref{eq:sketch} for the subspace \(\mathbb{R}^{n\times m}\) with high probability.
    Also, assume that for each restart index \(i \geq 1\) there exists \(k^{(i)} \leq m\) such that the computed preconditioned basis \(\hat Z_{1:k^{(i)}}^{(i)} \in \mathbb{R}^{n\times k^{(i)}}\) satisfies~\eqref{eq:compute-W}, \eqref{eq:compute-MLb}, \eqref{eq:compute-C}, \eqref{eq:compute-g}, \eqref{eq:sls}, \eqref{eq:updatex} with \(\epsW^{(i)}\), \(\epssls^{(i)}\), \(\epsapplyS^{(i)}\), \(\epsx^{(i)} \leq \bigO(\macheps)\),
    and that
    \begin{equation} \label{eq:thm:stability_restarted_sgmres:key_dimension}
    \begin{split}
        &\bigO(\macheps) \sqrt{\frac{1 - \epsilon}{1 + \epsilon}} \kappa\bigl(\bmat{\phi r^{(i)}& \hat W_{1:k^{(i)}}^{(i)}}\bigr) \geq 1,\qquad
        \bigO(\macheps)\sqrt{\frac{1 + \epsilon}{1 - \epsilon}} \kappa(\hat W_{1:k^{(i)}}^{(i)}) \leq \frac{1}{2}, \\
        &\text{and}\quad
        \bigO(\macheps) \kappa(\hat Z_{1:k^{(i)}}^{(i)}) \leq \frac{1}{2}.
    \end{split}
    \end{equation}
    
    If
    \begin{equation} \label{eq:thm:stability_restarted_sgmres:lambda1}
        \Lambda^{(i)}_{1} = \frac{1 + \epsilon}{1 - \epsilon} \bigO(\macheps)\kappa(M_L) \bigl(1 + \Bar{\tau}_{k^{(i)}}\bigr) \ll 1,
    \end{equation}
    holds for \(i\) with
    \[
        \Bar{\tau}_{k^{(i)}} = \frac{\norm{\hat Z^{(i)}_{1:k^{(i)}}} \normF{A} \norm{\hat y^{(k^{(i)})}}}{\norm{A\hat{d}^{(i)}}},
    \]
    then, with high probability, the backward error is reduced in the \(i\)-th restart cycle by a factor \(\Lambda^{(i)}_{1}\) until it holds that
    \begin{equation} \label{eq:restarted_back_err}
           \frac{\norm{b - A \hat x}}{\normF{A} \norm{\hat x} + \norm{b}} \leq \bigO(\macheps).
    \end{equation}
\end{theorem}

\begin{proof}
    Similar to the proof of~\cite[Theorem 4.1]{BHMV2024}, we will employ~\cite[Lemma 4.2]{BHMV2024}, which requires verifying the assumption~\cite[(4.14)]{BHMV2024}, i.e., bounding \(\norm{\hat{r}^{(i)} - A\hat{d}^{(i)}}\) by \(\norm{b - A\hat{x}^{(i)}}\) and \(\norm{b}+\normF{A}\norm{\hat{x}^{(i)}}\).
    From Lemma~\ref{lem:backward-sGMRES}, we have
    \begin{equation} \label{eq:thm:restart:norm_r-Ad}
        \begin{split}
            \norm{\hat{r}^{(i)} - A\hat{d}^{(i)}}
            &\leq \frac{1 + \epsilon}{1 - \epsilon} \bigO(\macheps)\kappa(M_L) \biggl(\frac{\norm{\hat Z^{(i)}_{1:k^{(i)}}} \normF{A} \norm{\hat y^{(k^{(i)})}}}{\norm{A\hat{d}^{(i)}}} \norm{A\hat{d}^{(i)}} + \norm{\hat{r}^{(i)}}\biggr) \\
            &\leq \frac{1 + \epsilon}{1 - \epsilon} \bigO(\macheps)\kappa(M_L) \biggl(\Bar{\tau}_{k^{(i)}} \bigl(\norm{\hat{r}^{(i)} - A\hat{d}^{(i)}} + \norm{\hat{r}^{(i)}}\bigr)
            + \norm{\hat{r}^{(i)}}\biggr) \\
            &\leq \frac{1 + \epsilon}{1 - \epsilon} \bigO(\macheps)\kappa(M_L) \biggl(\bigl(1 + \Bar{\tau}_{k^{(i)}}\bigr) \norm{\hat{r}^{(i)}}
            + \Bar{\tau}_{k^{(i)}} \norm{\hat{r}^{(i)} - A\hat{d}^{(i)}}\biggr).
        \end{split}
    \end{equation}
    By moving the term related to \(\norm{\hat{r}^{(i)} - A\hat{d}^{(i)}}\) from the right-hand side of~\eqref{eq:thm:restart:norm_r-Ad} to the left-hand side and using assumption~\eqref{eq:thm:stability_restarted_sgmres:lambda1}, we further obtain
    \begin{equation}
        \begin{split}
            \norm{\hat{r}^{(i)} - A\hat{d}^{(i)}}
            &\leq \frac{1 + \epsilon}{1 - \epsilon} \bigO(\macheps)\kappa(M_L) \bigl(1 + \Bar{\tau}_{k^{(i)}}\bigr) \norm{\hat{r}^{(i)}} \\
            &\leq \frac{1 + \epsilon}{1 - \epsilon} \bigO(\macheps)\kappa(M_L) \bigl(1 + \Bar{\tau}_{k^{(i)}}\bigr) \bigl(\norm{b - A\hat{x}^{(i)}} \\
            &\quad+ \bigO(\macheps)(\norm{b} + \normF{A}\norm{\hat{x}^{(i)}})\bigr),
        \end{split}
    \end{equation}
    which gives~\cite[(4.14)]{BHMV2024}.
    Together with~\cite[Lemma 4.2]{BHMV2024}, we can draw the conclusion.
\end{proof}

\begin{remark} \label{remark:thm-restart}
    When no preconditioning is used, i.e., \(M_L = M_R = I\) and \(\hat Z_{1:k^{(i)}}^{(i)} = \hat B_{1:k^{(i)}}^{(i)}\), if it holds for the key dimension \(k^{(i)}\) that
    \begin{equation} \label{eq:thm:stability_restarted_sgmres_nopre:lambda1}
        \Lambda^{(i)}_{1} = \bigO(\macheps) \frac{1 + \epsilon}{1 - \epsilon} \biggl(1+\frac{\norm{\hat B^{(i)}_{1:k^{(i)}}} \normF{A} \norm{\hat y^{(k^{(i)})}}}{\norm{A\hat{d}^{(i)}}}\biggr) \ll 1,
    \end{equation}
    then, with high probability, the backward error is reduced in the \(i\)-th restart cycle by a factor \(\Lambda^{(i)}_{1}\) until it holds that
    \begin{equation*}
        \frac{\norm{b - A \hat x}}{\normF{A} \norm{\hat x} + \norm{b}} \leq \bigO(\macheps).
    \end{equation*}
\end{remark}

\begin{remark} \label{remark:thm:stability_restarted_sgmres}
As a consequence of Theorem~\ref{thm:stability_restarted_sgmres}, it is clear that a restarted variant of sketched GMRES (Algorithm \ref{alg:restarted_sGMRES}) may allow us to recover the stability of sGMRES when \(\Bar{\tau}_{k^{(i)}}\) is moderate.
However, it should be noted that the assumptions, particularly~\eqref{eq:thm:stability_restarted_sgmres:key_dimension}, of this theorem are overly pessimistic.
In general, selecting a sufficiently large \(m\) is required to guarantee~\eqref{eq:thm:stability_restarted_sgmres:key_dimension}, but, in practice, we typically select a relatively small value for the maximum iterations per restart cycle, which makes it difficult to satisfy~\eqref{eq:thm:stability_restarted_sgmres:key_dimension}.
\end{remark}

As discussed in Remark~\ref{remark:thm:stability_restarted_sgmres}, for restarted sGMRES, there are two possible reasons for instability: One is due to choosing an \(m\) that is too small; the other is due to using a value of \(t\) that is too small, which has an influence on \(\Bar{\tau}_{k^{(i)}}\).
Since the choice of \(m\) depends mainly on storage in practice, here we only consider recovering from the instability caused by a small \(t\).
Employing a small \(t\) can result in considerable acceleration, though it might also introduce instability issues.

Next, our aim is to propose an adaptive strategy to recover from the instability of sGMRES caused by using a small \(t\).
The definition of \(\Lambda_1^{(1)}\) in Theorem~\ref{thm:stability_restarted_sgmres} relying on \(\Bar{\tau}_{k^{(i)}}\), which influences the reduced factor of the backward error in the \(i\)-th restart cycle, implies that
care should be taken to ensure \(\Bar{\tau}_{k^{(i)}}\) does not get too large as the iteration progresses.
In fact, as noted in \cite{NT2024}, due to the subspace embedding property \eqref{eq:sketch}, \(\Bar{\tau}_{i}\) can be inexpensively estimated as 
\begin{equation*}
   \tilde{\tau}_{i} := \frac{\norm{ S Z_{1:i}} \normF{A} \norm{y^{(i)}}}{\norm{SAZ_{1:i}y^{(i)}}}.
\end{equation*}
The numerical stability of sGMRES can thus be controlled by monitoring \(\tilde{\tau}_{i}\) at each iteration, and doubling the value of \(t\) when it grows above some specified tolerance.
We then proceed to double \(t\) again only if \(\tilde{\tau}_{i}\) is still large and \(\tilde{\tau}_{i}\) significantly exceeds \(\tilde{\tau}_{i-1}\), indicating that \(t\) remains insufficiently large.
Restarted sGMRES combined with this adaptive strategy is presented in Algorithm~\ref{alg:sGMRES_adapt}.

\begin{algorithm}[!tb]
\begin{algorithmic}[1]
    \caption{One cycle of preconditioned sketched GMRES (sGMRES) \label{alg:sGMRES_adapt} with adaptive truncation size}
    \REQUIRE
    A matrix \(A \in \mathbb R^{n\times n}\), a right-hand side \(b \in \mathbb R^{n}\), an initial approximation \(x^{(0)} \in \mathbb R^{n}\), the maximal number of iterations \(m\), a subspace embedding \(S \in \mathbb R^{s\times n}\) (\(s>m\)), left and right preconditioners \(M_L\) and \(M_R\), and truncation parameter $t$, \(\texttt{tol}_{\tau}\) and \(\texttt{tol}_{B}\) for the restarting strategy.
    \ENSURE
    A computed solution \(x \in \mathbb R^{n}\) approximating the solution of \(A x = b\).
    \STATE \(r^{(0)}\gets b - Ax^{(0)}\) and \(g \gets S M_L^{-1} r^{(0)}\). 
    \STATE \(\tilde x^{(0)} \gets Sx^{(0)}\).
    \STATE \(B_{1} \gets r^{(0)}/\norm{r^{(0)}}\).
    \FOR{\(i = 1:m\)}
        \STATE Generate the \(i\)-th column \(B_{i}\) of the basis by the Arnoldi process, and \(Z_{i}\gets M_R^{-1} B_i\).
        \STATE \(W_{i} \gets M_L^{-1} A Z_{i}\). 
        \STATE Sketch reduced matrix \(C_{1:i} \gets S W_{1:i}\). 
        \STATE Compute the QR factorization \(C_{1:i} = U_{1:i} T_{1:i}\) using \(C_{1:i-1} = U_{1:i-1} T_{1:i-1}\), where \(U_{1:i}\) is an orthonormal matrix and \(T_{1:i}\) is an upper triangular matrix.
        \STATE Solve the triangular system \(T_{1:i} y^{(i)} = U_{1:i}\trans g\) to obtain \(y^{(i)}\).
        \STATE \(\tilde Z_i = SZ_i\) and \(\tilde{d}^{(i)} = C_{1:i} y^{(i)}\).
        \IF{the stopping criterion is satisfied}
            \RETURN \(x = x^{(i)} \gets x^{(0)} + d^{(i)}\) with \(d^{(i)} = Z_{1:i} y^{(i)}\).
        \ENDIF
        \STATE $\tilde{\tau}_{i} \leftarrow \norm{\tilde Z_{1:i}} \normF{A}\norm{y^{(i)}}/\norm{\tilde{d}^{(i)}}$.
        \IF{\(\texttt{tol}_{\tau}\cdot \tilde{\tau}_{i} \geq 1\)}
            \IF{\(\tilde{\tau}_{i} > 1.1\cdot \tilde{\tau}_{i-1}\)}
                \STATE{$t \leftarrow \min\{i+1, 2 t\}$}.
            \ENDIF
        \ENDIF
    \ENDFOR
\end{algorithmic}
\end{algorithm}

\section{Numerical Experiments}
\label{sec:experiments}
In this section, we present numerical experiments that demonstrate the results of our stability analysis.
All experiments were performed in MATLAB R2023a on a Windows~11 HP laptop with an 11th Gen Intel(R) Core processor with 2.80\,GHz and 8\,GB of RAM. We used (with modifications) the open source code provided by G{\"u}ttel and Simunec as part of their recent work \cite{guttel2024sketch}.
Our code is available at https://github.com/burkel8/stability\_sGMRES.
The sketching operator used in all experiments is the subsampled random Hadamard transform (see, e.g., ~\cite{balabanov2025randomized}).

We consider solving a linear system $A x = b$ with sGMRES using both the truncated Arnoldi process and the sketch-and-select Arnoldi method (sGMRES-trunc and sGMRES-ssa, respectively, in the plots).
Although a few variants of the sketch-and-select Arnoldi algorithm have been proposed in \cite{guttel2024sketch}, we use the 
\texttt{sGMRES-ssa-pinv} method, as it was shown to produce the best conditioned basis in their experiments. We also ran the experiments with \texttt{sGMRES-ssa-greedy} and obtained similar results, and thus do not include this method here. We note that the different variations of sketch-and-select Arnoldi differ in how the $t$ vectors with which to orthogonalize against are chosen, and we refer to \cite{guttel2024sketch} for more details on these methods. In these experiments, we always choose the same value of $t$ in the sketch-and-select Arnoldi algorithm that we use as the truncation parameter in the truncated Arnoldi method. All implementations of standard GMRES use the modified Gram-Schmidt Arnoldi process. We do not use any preconditioner, i.e., \(M_L=M_R=I\), in the experiments except for Figure~\ref{fig:adapt_sherman2_pre}. Finally, we drop the $\hat{\cdot}$ notation in this section, as all quantities are computed quantities.

\subsection{Illustration of the backward error bound}
In this subsection, we use three different examples to illustrate the results of our analysis shown in Sections~\ref{sec:stability} and~\ref{sec:restart-sGMRES}.
For all sketching techniques we take $s = 2(m+1)$, where $m$ is the maximum length of each Arnoldi cycle, as suggested in \cite{NT2024}. We display the backward error, $\tau_{i}$ from \eqref{eq:lem-res}, and the condition numbers of $B_{1:i}$ and $A B_{1:i}$ (approximated using the condition numbers of $S B_{1:i}$ and $SAB_{1:i}$, respectively) at each iteration.
These condition number approximations are denoted as $\kappa(B_{1:i})$ and $\kappa(A B_{1:i})$ in the plots.

\subsubsection{Attainable worst case for sGMRES}
In the first example, the matrix $A$ is a $400 \times 400$ random banded matrix with condition number $\kappa(A) = 10$, obtained using the MATLAB \texttt{gallery} function with \texttt{randsvd} as input.
The right-hand side is taken to be the third singular vector of $A$.
We run $m = 400$ iterations of standard GMRES and sGMRES with both truncated Arnoldi and sketch-and-select Arnoldi with $t = 2$, and do not use any randomized subspace embedding for this example, i.e., \(S\) is chosen to be a \(400\)-by-\(400\) identity matrix.
Theoretically, the identity matrix is the most stable choice for \(S\), but is not an appropriate choice of sketching operator in practice. Here, we use the identity matrix solely to illustrate the worst case for sGMRES, i.e., the case where backward error of sGMRES barely decreases due to the ill-conditioned basis $B_{1:i}$.

\begin{figure}[!tb]
    \centering
    \includegraphics[width= \textwidth]{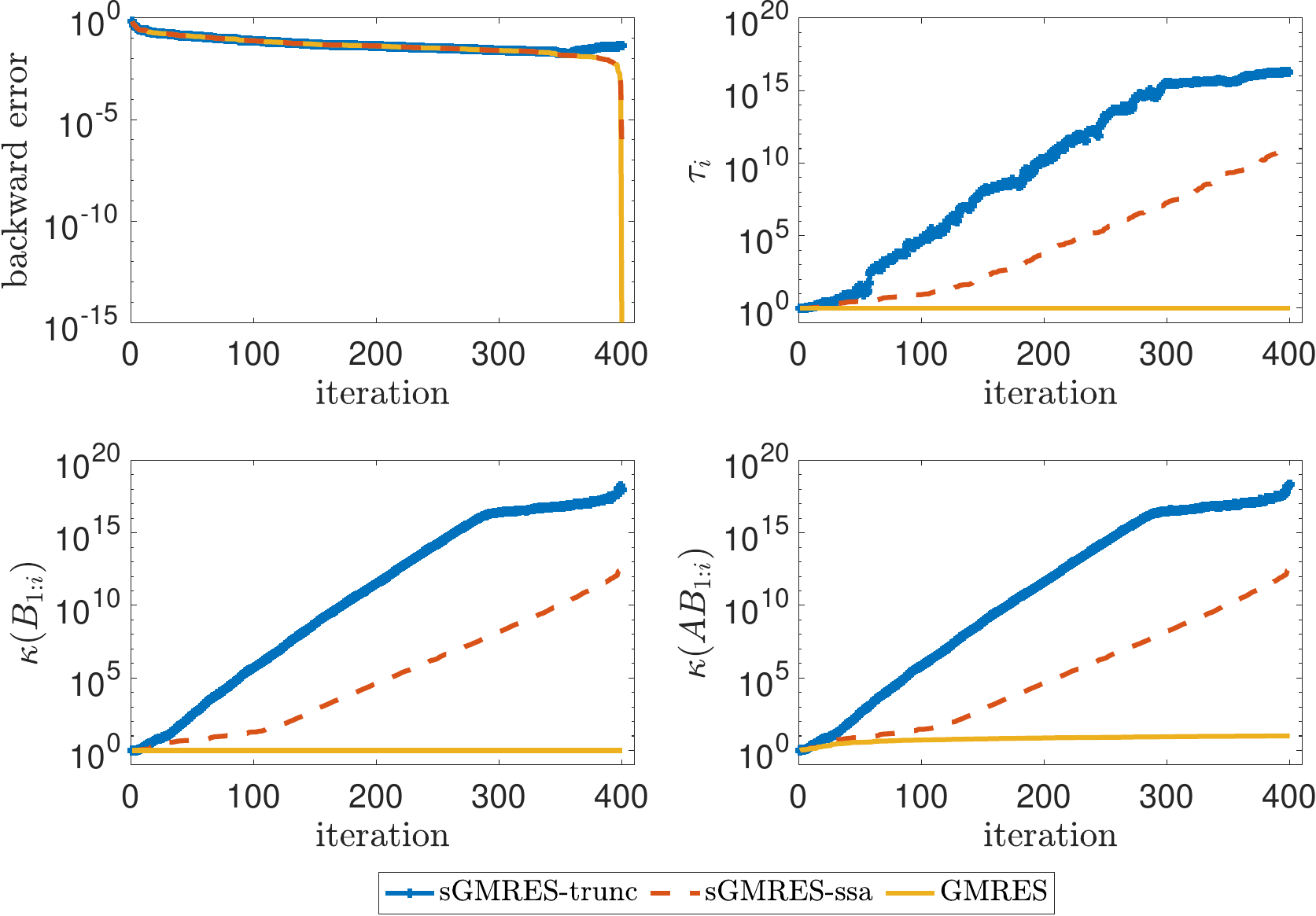}
    \caption{Example 1: From the left to right, the plots show the backward error,  $\tau_{i}$ and the approximate condition numbers of basis $B_{i}$ and $A B_{i}$, monitored every iteration, for different implementations of sGMRES, and an implementation of standard GMRES.}
    \label{fig:exp1}
\end{figure}

The results are shown in Figure \ref{fig:exp1}, where we see that the condition number of $B_{1:i}$ and $\tau_{i}$ generated by the truncated Arnoldi grows to $10^{16}$, leading to a poor backward error result for \texttt{sGMRES-trunc}.
In the case of \texttt{sGMRES-ssa}, the backward error is approximately \(10^{-7}\) as both the condition number of $B_{1:i}$ and $\tau_{i}$ are not too large.
For GMRES, it can be seen that the condition number of $B_{1:i}$ remains small, as does $\tau_{i}$, yielding a very small backward error after $400$ iterations. 
Thus, we can see that the bound in Theorem~\ref{thm:stability_sgmres} can be attainable for some cases.

\subsubsection{Comparison between the bounds shown in Theorem~\ref{thm:stability_sgmres} and Lemma~\ref{lem:backward-sGMRES}}
In the second example, the matrix $A$ is the \texttt{Norris/torso3} matrix of size \(n = 259,156\), and the right-hand side vector is randomly generated according to the standard normal distribution. 
We take $m=300$ and $t = 2$ for this example.
We use this example to demonstrate that the bound outlined in Lemma~\ref{lem:backward-sGMRES} is sharper than the one mentioned in Theorem~\ref{thm:stability_sgmres}, and proves useful in certain situations.

\begin{figure}[!tb]
    \centering
    \includegraphics[width= \textwidth]{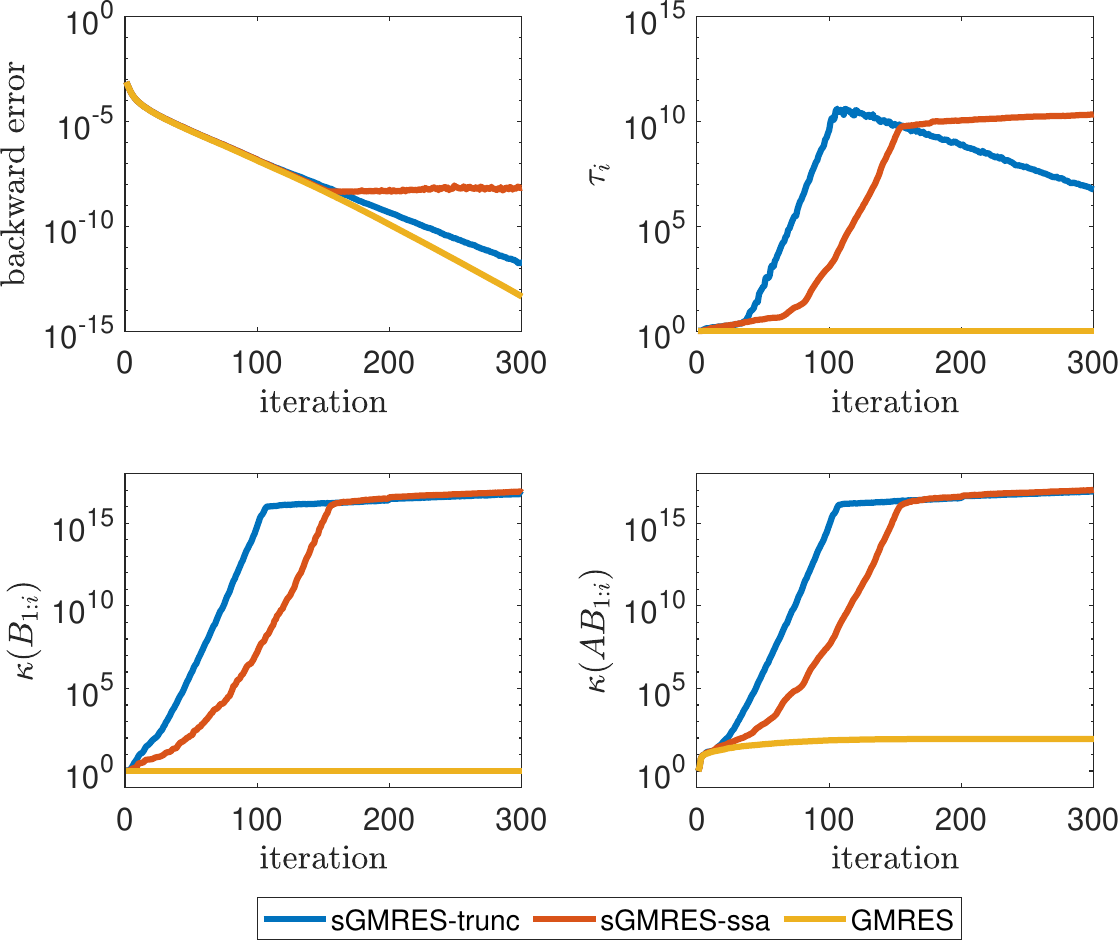}
    \caption{Example 2: From the left to right, the plots show the backward error,  $\tau_{i}$ and the approximate condition numbers of the $B_{i}$ and $A B_{i}$, monitored every iteration, for different implementations of sGMRES, and an implementation of standard GMRES.}
    \label{fig:exp2}
\end{figure}

As shown in Figure~\ref{fig:exp2}, the behavior of the backward error for both \texttt{sGMRES-trunc} and \texttt{sGMRES-ssa} diverges; notably, while the backward error for \texttt{sGMRES-ssa} stops improving after about 160 iterations, \texttt{sGMRES-trunc} continues to reduce its error reaching a small level, despite both methods having ill-conditioned $B_{1:i}$.
This indicates that for this scenario, the bounds shown by Theorem~\ref{thm:stability_sgmres}, specifically~\eqref{eq:res-xk} and~\eqref{eq:remark1_bound}, are overly pessimistic and fail to account for the discrepancy observed. 
We further notice that $\tau_{i}$ for \texttt{sGMRES-trunc} is considerably smaller than for \texttt{sGMRES-ssa}, explaining the superior backward error performance of \texttt{sGMRES-trunc}, as anticipated by Lemma~\ref{lem:backward-sGMRES}.
We also observe but do not find a clear theoretical explanation for the phenomenon that \(\tau_i\) can even decrease after \(\kappa(AB)\) and \(\kappa(B)\) are both very ill-conditioned, which means that solving the least squares problem is totally unstable.

\subsubsection{Comparison between different sGMRES variants with and without restarting}
In the third example, the matrix $A$ is the \texttt{Norris/stomach} matrix of size $n = 213,360$ obtained from the SuiteSparse matrix collection \cite{Davis2011} (also tested in~\cite{guttel2024sketch}). The right-hand side vector is randomly generated according to the standard normal distribution.
We run $\texttt{nrestarts}=5$ cycles of restarted sGMRES with the truncated Arnoldi process, and sketch-and-select Arnoldi, with $m = 150$ iterations per cycle, totaling $750$ iterations.
We also run restarted GMRES with the same parameters.
Additionally, we perform sGMRES with truncated Arnoldi and sketch-and-select Arnoldi with $m = 750$ iterations (without restarts). All implementations of sGMRES use a truncation parameter $t = 3$.

\begin{figure}[!tb]
    \centering
\includegraphics[width= \textwidth]{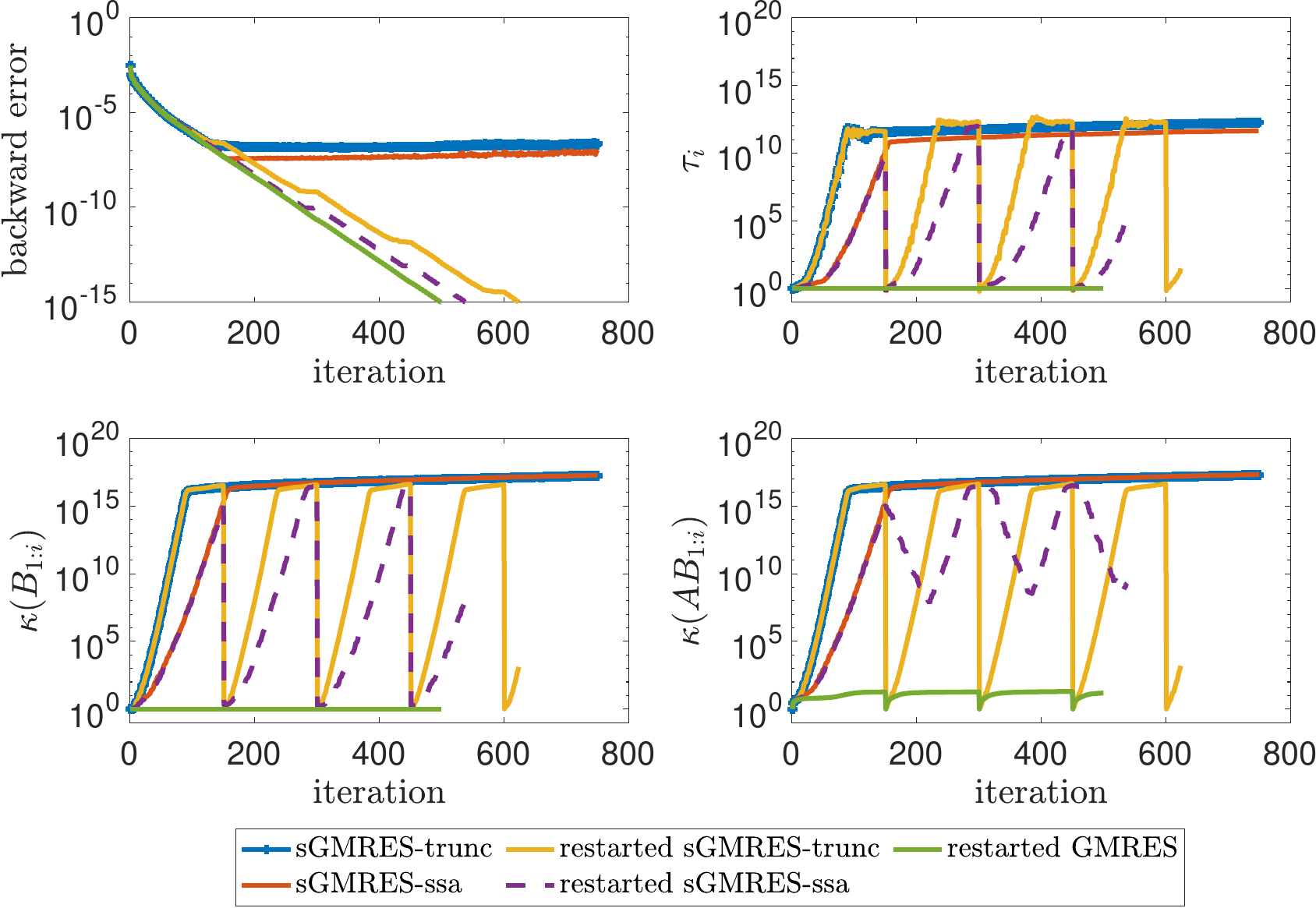}
    \caption{Example 3: From the left to right, the plots show the backward error, $\tau_{i}$, and the approximate condition numbers of $B_{i}$ and $A B_{i}$, monitored every iteration, for different implementations of sGMRES, and an implementation of standard GMRES.}
    \label{fig:exp3}
\end{figure}

The results are shown in Figure \ref{fig:exp3}, where we observe that for sGMRES with truncated Arnoldi and sketch-and-select Arnoldi (without restarts), \(\tau_i\) and the condition numbers of both $B_{1:i}$ and $A B_{1:i}$ grow to about $10^{16}$, causing the backward error to stagnate, in accordance with the results of Theorem~\ref{thm:stability_sgmres} and Lemma~\ref{lem:backward-sGMRES}.
The restarted implementation helps regain stability, leading to a continued reduction in backward error with each restart, eventually reaching \(\bigO(\macheps)\).
This is because restarting can remove the influence of the condition number of $B_{1:i}$ as predicted in Theorem~\ref{thm:stability_restarted_sgmres}.

In Figure \ref{fig:exp3}, we also see that the backward error for sGMRES restarted with a sketch-and-select Arnold reduces faster than for sGMRES restarted with truncated Arnoldi in each cycle.
Although the condition numbers of $B_{1:i}$ produced by these two methods eventually become large, the values of $\tau_{i}$ are small.
Thus, they can produce small backward errors, which again confirms that the bound in Lemma~\ref{lem:backward-sGMRES} is sharper than in Theorem~\ref{thm:stability_sgmres}.

\subsection{sGMRES with adaptive truncation size}
In the final experiment, we test the adaptive strategy discussed in Section~\ref{sec:sketched_GMRES}. In particular, we test
\begin{itemize}
    \item \textbf{restarted GMRES:} restart every \(m\) iterations.
    \item \textbf{restarted sGMRES:} restart every \(m\) iterations with fixed parameter \(t = 1\).
    \item \textbf{restarted sGMRES with adapt \(t\):} restart every \(m\) iterations with adaptive parameter \(t\) beginning with \(t=1\) and \(\texttt{tol}_{\tau} = \macheps\), i.e., Algorithm~\ref{alg:sGMRES_adapt}.
    \item \textbf{restarted GMRES with \(m = t_{\max}\):} restart every \(t_{\max}\) iterations, where \(t_{\max}\) denotes the maximum truncation size employed in the ``restarted sGMRES with adapt \(t\)'' method.
\end{itemize}

\begin{figure}[!tb]
    \centering
    \includegraphics[width= \textwidth]{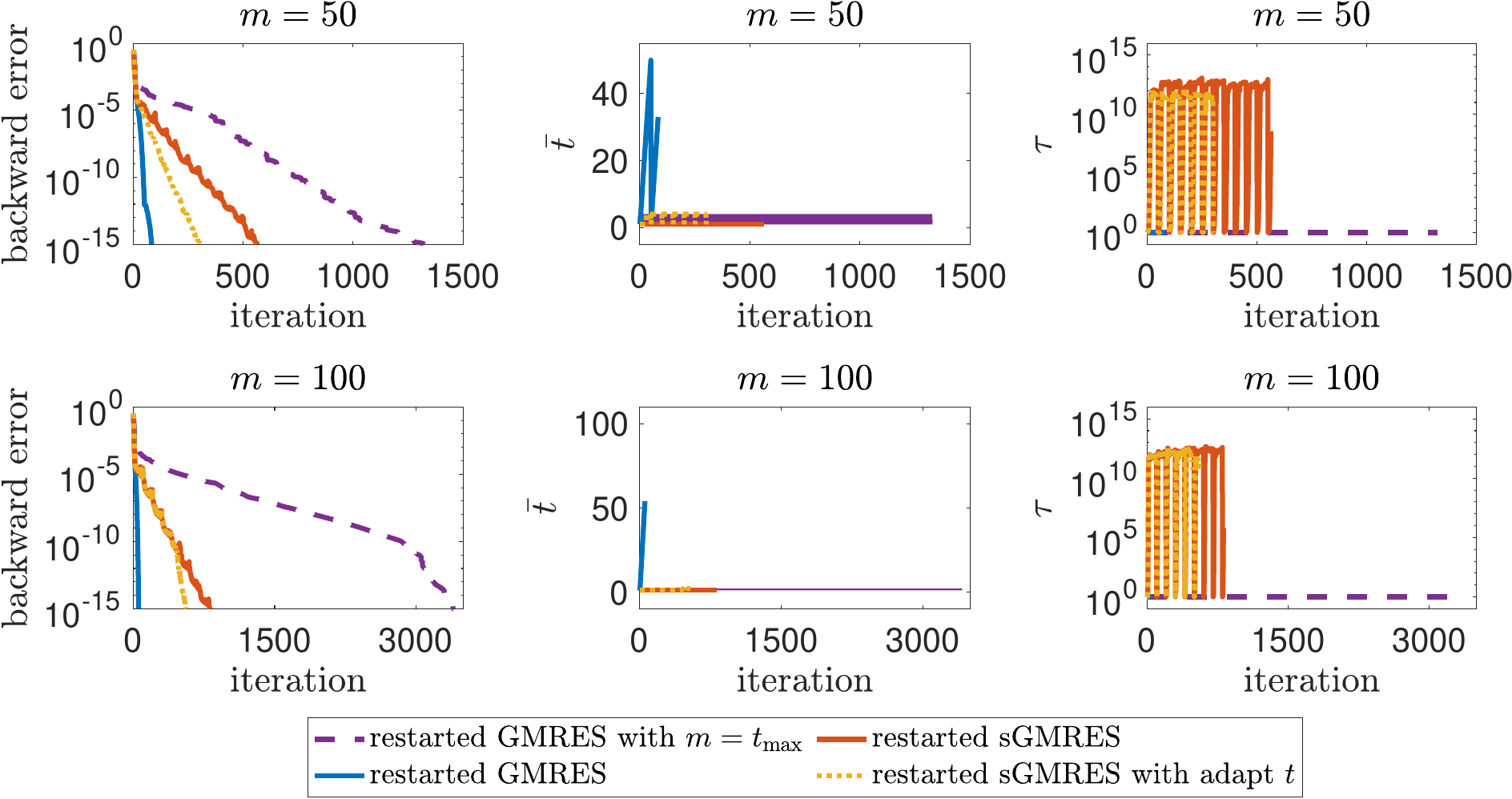}
    \caption{Adaptive truncation size experiment. From left to right, we plot the backward error, the value of \(\bar{t}\), and the value of \(\tau\) for \texttt{fs\_760\_1}, where \(\bar{t} = \max(t, i)\) represents the number of columns against which \(B_i\) needs to perform orthogonalization in Line~\ref{line:algo_trunc:forloop} of Algorithm~\ref{alg:trunc_Arnoldi} and \(t_{\max}\) denotes the maximum truncation size employed in the ``restarted sGMRES with adapt \(t\)'' method.}
    \label{fig:adapt_fs760}
\end{figure}

\begin{figure}[!tb]
    \centering
    \includegraphics[width= \textwidth]{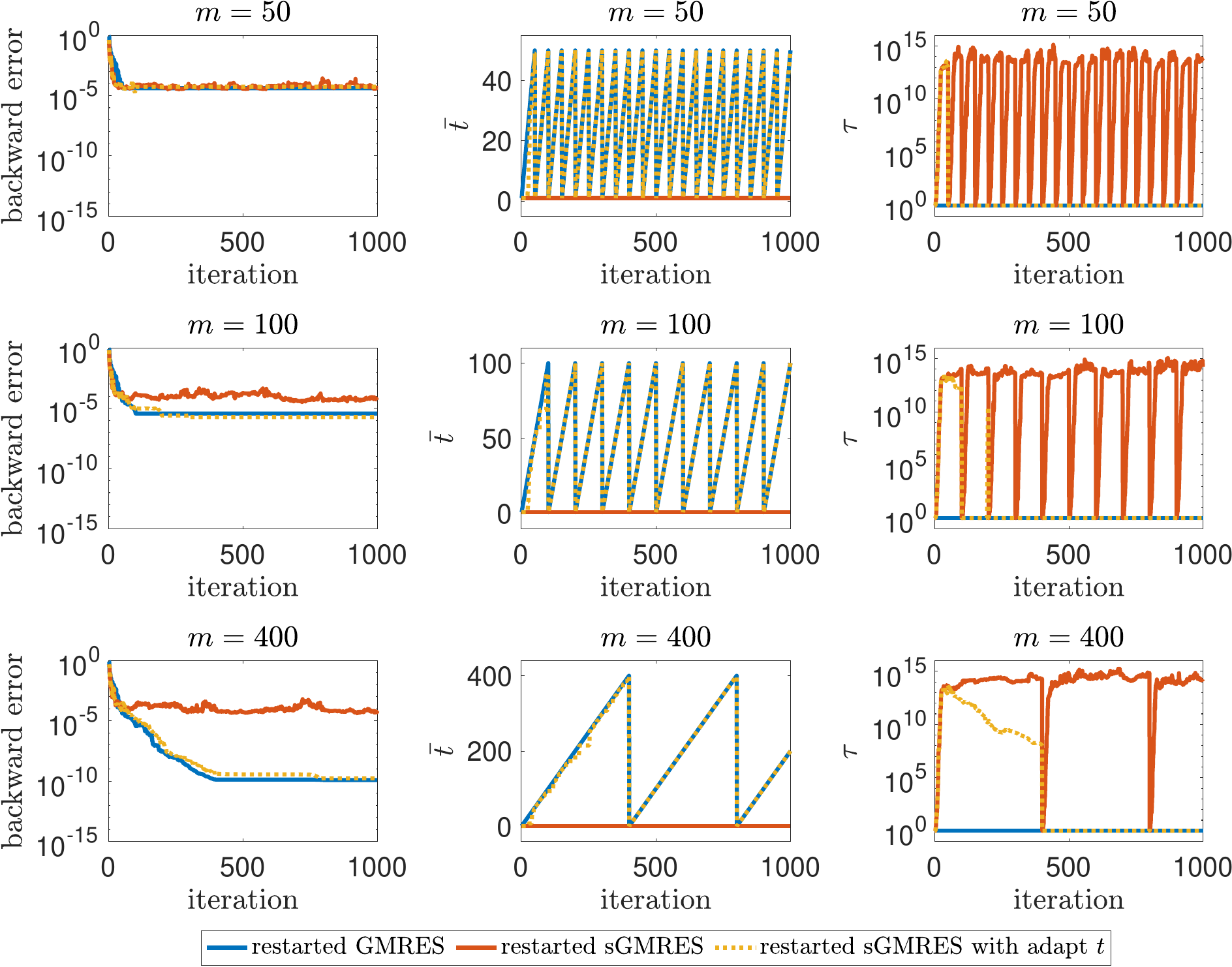}
    \caption{Adaptive truncation size experiment. From left to right, we plot the backward error, the value of \(\bar{t}\), and the value of \(\tau\) for \texttt{sherman2}, where \(\bar{t} = \max(t, i)\) represents the number of columns against which \(B_i\) needs to perform orthogonalization in Line~\ref{line:algo_trunc:forloop} of Algorithm~\ref{alg:trunc_Arnoldi}.}
    \label{fig:adapt_sherman2}
\end{figure}

Figures~\ref{fig:adapt_fs760} and~\ref{fig:adapt_sherman2} illustrate tests of these algorithms on matrices \texttt{fs\_760\_1} and \texttt{sherman2}, whose condition numbers are \(5.49\times 10^3\) and \(9.64\times 10^{11}\), respectively. 
In each example, the right-hand side is generated randomly following the standard normal distribution, with \(m\) values of \(50\) and \(100\), and \(\texttt{nrestarts} = 1.5\cdot 10^3/m\).
It should be noted that restarted GMRES with \(m = t_{\max}\) is examined primarily to enable a direct comparison between restarted sGMRES with adaptive \(t\) and restarted
GMRES with similarly low orthogonalization costs.
In Figure~\ref{fig:adapt_sherman2}, we omit the result for ``restarted GMRES with \(m = t_{\max}\)'' because it coincides with that of ``restarted GMRES''.

In Figure~\ref{fig:adapt_fs760}, both restarted GMRES and sGMRES achieve small backward errors, with restarted sGMRES using the adaptive strategy accelerating convergence compared to a fixed \(t\) approach by using a slightly larger \(t\).
As shown in the middle subfigure of Figure~\ref{fig:adapt_fs760}, it can be expected that the performance of restarted sGMRES, using adaptive truncation size, tends towards the performance of restarted sGMRES.
To further illustrate the potential performance advantage of our adaptive strategy, we also compare the number of floating-point operations required for orthogonalization in Table~\ref{table:flops}.
From the table, we observe that restarted sGMRES, both with and without the adaptive strategy, entails a comparable but substantially lower orthogonalization cost than restarted GMRES.

We next present a difficult example for restarted sGMRES in Figure~\ref{fig:adapt_sherman2}.
In this situation, restarted sGMRES struggles to attain a small backward error.
As indicated in Theorem~\ref{thm:stability_restarted_sgmres}, the attainable backward error is affected by two conditions, namely \eqref{eq:thm:stability_restarted_sgmres:key_dimension} and~\eqref{eq:thm:stability_restarted_sgmres:lambda1}.
For \(m = 50\) in Figure~\ref{fig:adapt_sherman2}, the value of \(m\) is too small for \eqref{eq:thm:stability_restarted_sgmres:key_dimension} to be satisfied.
Consequently, even if we use restarted sGMRES with the adaptive strategy or restarted GMRES, both of which ensure a small \(\tau\) so that~\eqref{eq:thm:stability_restarted_sgmres:lambda1} holds, restarting still fails to improve the attainable backward error.
In contrast, when \(m = 100\) and \(m = 400\), the larger values of \(m\) allow~\eqref{eq:thm:stability_restarted_sgmres:key_dimension} to be satisfied.
In this situation, \(\tau\), associated with~\eqref{eq:thm:stability_restarted_sgmres:lambda1}, becomes the dominant factor determining the attainable backward error.
Figure~\ref{fig:adapt_sherman2} clearly shows that, for \(m = 100\), our adaptive strategy can enhance the poor backward error that arises from a small \(t\) by adaptively increasing \(t\) to reduce \(\tau\) such that~\eqref{eq:thm:stability_restarted_sgmres:lambda1} is satisfied.
This effect is even more pronounced for \(m = 400\).
In this case, \(t\) is repeatedly doubled up to \(m\) during the early stage, thus restarted sGMRES with the adaptive strategy utilizes a comparable orthogonalization cost as standard GMRES.
By comparing the experiments shown in Figures~\ref{fig:adapt_fs760} and~\ref{fig:adapt_sherman2}, we observe that restarted GMRES also requires a large \(m\) in order to obtain a small backward error. 
According to Theorem~\ref{thm:stability_restarted_sgmres}, \(\bigO(\macheps)\sqrt{\frac{1 + \epsilon}{1 - \epsilon}} \kappa(\hat W_{1:k^{(i)}}^{(i)}) \leq 1/2\) in~\eqref{eq:thm:stability_restarted_sgmres:key_dimension} depends on \(\kappa(\hat W_{1:k^{(i)}}^{(i)})\), whose upper bound is determined by \(\kappa(A)\) and \(\kappa(\hat Z_{1:k^{(i)}}^{(i)})\).
This implies that, in the worst case, it may be necessary to require \(\bigO(\macheps)\kappa(\hat Z_{1:k^{(i)}}^{(i)})\leq 1/(2\kappa(A))\) in order to guarantee~\eqref{eq:thm:stability_restarted_sgmres:key_dimension}.
Consequently, for well-conditioned matrices such as \texttt{fs\_760\_1}, \(\kappa(\hat Z_{1:k^{(i)}}^{(i)})\) can be relatively large while still satisfying \(\bigO(\macheps)\kappa(\hat Z_{1:k^{(i)}}^{(i)})\leq 1/(2\kappa(A))\), meaning that a relatively small truncation size usually suffices.
In contrast, for ill-conditioned matrices, a significantly larger small truncation size may be required.

\begin{figure}[!tb]
    \centering
    \includegraphics[width= \textwidth]{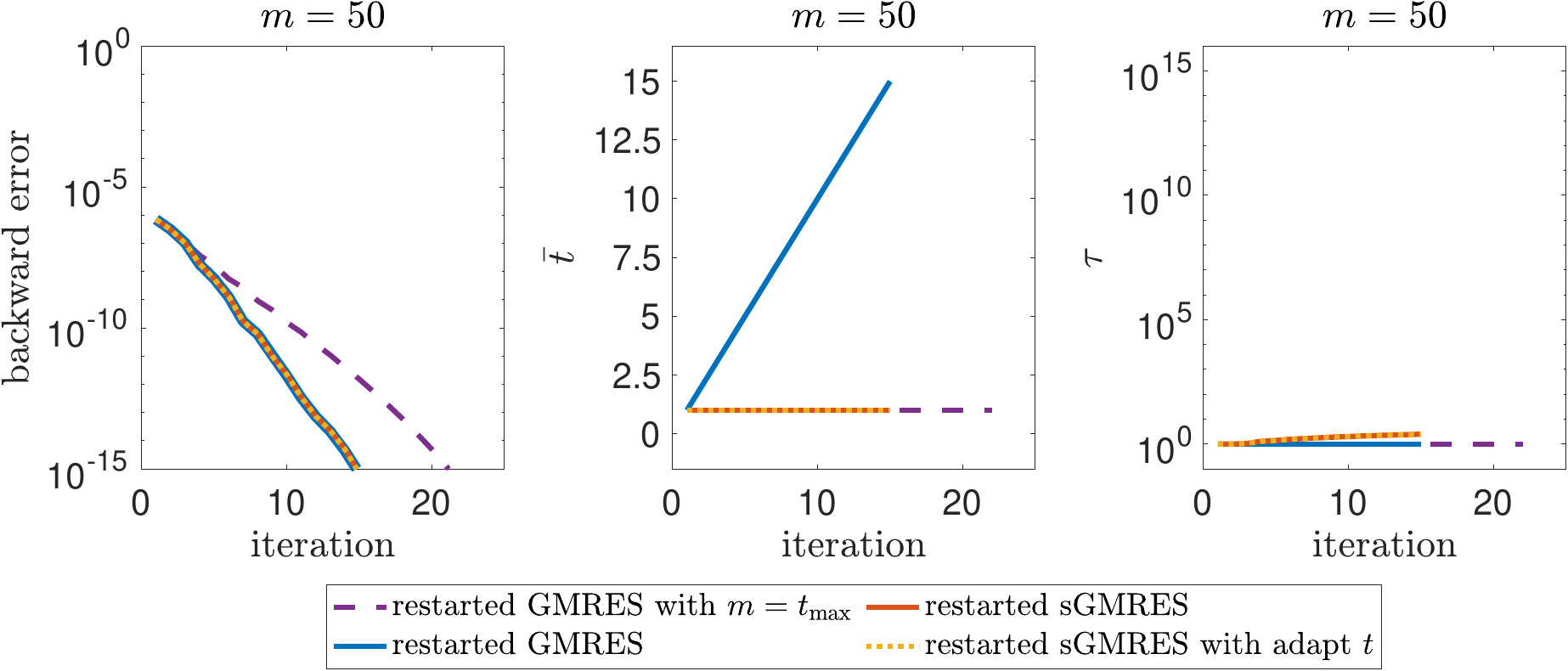}
    \caption{Adaptive truncation size experiment with preconditioning. From left to right, we plot the backward error, the value of \(\bar{t}\), and the value of \(\tau\) for \texttt{sherman2}, where \(\bar{t} = \max(t, i)\) represents the number of columns against which \(B_i\) needs to perform orthogonalization in Line~\ref{line:algo_trunc:forloop} of Algorithm~\ref{alg:trunc_Arnoldi} and \(t_{\max}\) denotes the maximum truncation size employed in the ``restarted sGMRES with adapt \(t\)'' method.}
    \label{fig:adapt_sherman2_pre}
\end{figure}

Furthermore, in Figure~\ref{fig:adapt_sherman2_pre}, we combine the preconditioners generated by the MATLAB command \(\texttt{ilu}\) with the adaptive strategy to test \texttt{sherman2}, i.e., the same matrix tested in Figure~\ref{fig:adapt_sherman2}.
Due to the efficiency of the preconditioners, a small \(m\), that is, \(m = 50\), is sufficient to achieve convergence.
From Figure~\ref{fig:adapt_sherman2_pre} together with Table~\ref{table:flops}, it can be seen that when preconditioning is used, restarted sGMRES with and without the adaptive strategy have similar behavior and require fewer orthogonalizations than standard GMRES.

In general, sGMRES with adaptive truncation size often performs comparably to standard GMRES regarding the accuracy, and has the potential to result in substantial time savings especially for the well-conditioned case.

\begin{table}[!tb]
\caption{Comparison of the relative number of floating-point operations used in orthogonalization.
In the table, we report the ratios of the number of floating-point operations required for the orthogonalization of ``restarted sGMRES,'' ``restarted sGMRES with adapt \(t\),'' and ``restarted GMRES with \(m = t_{\max}\)'' relative to that of ``restarted GMRES,'' respectively.}
\label{table:flops}
\begin{tabular}{c|c|cccc}
\hline
 & m & \begin{tabular}[c]{@{}c@{}}restarted\\  GMRES\end{tabular} & \begin{tabular}[c]{@{}c@{}}restarted\\  sGMRES\end{tabular} & \begin{tabular}[c]{@{}c@{}}restarted\\ sGMRES\\ with adapt \(t\)\end{tabular} & \begin{tabular}[c]{@{}c@{}}restarted\\ GMRES\\ with \(m = t_{\max}\)\end{tabular} \\ \hline
\multirow{2}{*}{\begin{tabular}[c]{@{}c@{}}\texttt{fs\_760\_1} \\ without pre\end{tabular}} & 50 & 1 & 0.3072 & 0.5959 & 1.7979 \\ \cline{2-6} 
 & 100 & 1 & 0.5468 & 0.4485 & 3.4323 \\ \hline
\begin{tabular}[c]{@{}c@{}}\texttt{fs\_760\_1} \\ with pre\end{tabular} & 50 & 1 & 0.1250 & 0.1250 & 0.1883 \\ \hline
\end{tabular}
\end{table}

\section{Conclusions}
\label{sec:conclusions}
In this paper, we performed a backward stability analysis of sketched GMRES and showed how the backward error is largely influenced by the condition number of the generated Krylov basis.
However, in~\cite[Figure 9]{guttel2024sketch}, it was demonstrated that the backward error of sGMRES can still decrease even when the basis becomes ill-conditioned.
To provide more insight into this finding from~\cite{guttel2024sketch}, we proposed a sharper bound that relies on \(\norm{\hat B_{1:k}} \norm{\hat y^{(k)}}/\norm{\hat x^{(k)}}\) rather than \(\kappa(\hat B_{1:k})\).
Nevertheless, it remains unclear why \(\norm{\hat B_{1:k}} \norm{\hat y^{(k)}}/\norm{\hat x^{(k)}}\) itself can even decrease—leading to a reduction in the backward error of sGMRES—when the basis becomes ill-conditioned.
We then showed how restarting can be used to stabilize a sketched GMRES implementation, and proposed an adaptive strategy to choose proper truncation sizes for sGMRES based on our analysis.
Our numerical experiments indicate that the convergence of sketched GMRES with adaptive truncation size often performs comparably to standard GMRES and has the potential to result in substantial time savings.
\section*{Acknowledgments}
All authors are supported by the Charles University Research Centre program No. UNCE/\-24/\-SCI/005 and by the European Union (ERC, inEXASCALE, 101075632). Views
and opinions expressed are those of the authors only and do not necessarily reflect those of the European Union or the European Research Council. Neither the European Union nor the granting authority can be held responsible for them. L.~Burke acknowledges
funding from the InterFLOP (ANR-20-CE46-0009), NumPEx
ExaMA (ANR-22-EXNU-0002), MixHPC (ANR-23-CE46-0005-01), and FPT-4 (ANR-24-CE46-
7572) projects of the French National Agency for Research (ANR).

\bibliographystyle{abbrvurl}
\bibliography{mybib}


\end{document}